\documentclass{amsart}
\usepackage{amssymb, amsmath}
\newtheorem{theorem}{Theorem}
\newtheorem{proposition}[theorem]{Proposition}
\newtheorem{question}[theorem]{Question}
\newtheorem{problem}[theorem]{Problem}
\newtheorem{example}[theorem]{Example}
\newtheorem{definition}[theorem]{Definition}
\newtheorem{lemma}[theorem]{Lemma}
\newtheorem{corollary}[theorem]{Corollary}
\newtheorem{conjecture}[theorem]{Conjecture}

\title[Variations of selective separability II]{Variations of selective separability II:\\
discrete sets and the influence of convergence and maximality}
\author{Angelo Bella}

\address{Dipartimento di matematica\\
Citt\'a Universitaria\\ Viale A. Doria 6 \\ 98125 Catania\\ Italy}

\email{bella@dmi.unict.it}

\author{Mikhail Matveev}

\address{Department of mathematical Sciences\\ George Mason University\\
4400 University Drive\\ Fairfax, VA\\ 22030, USA}

\email{mmatveev@gmu.edu}

\author{Santi Spadaro}

\address{Department of Mathematics\\ Ben Gurion University of the Negev\\
Be'er Sheva \\ 84105 Israel}

\email{santi@cs.bgu.ac.il, santispadaro@yahoo.com}

\thanks{The third author was supported by the Center for Advanced Studies in Mathematics at Ben Gurion University}



\date{}
\begin{document}
\begin{abstract}
A space $X$ is called \emph{selectively separable} (\emph{R-separable}) if for every
sequence of dense subspaces $(D_n : n\in\omega)$ one can pick finite
(respectively, one-point) subsets $F_n\subset D_n$ such that
$\bigcup_{n\in\omega}F_n$ is dense in $X$. These properties are much stronger than separability, but are equivalent to it in the presence of certain convergence properties. For example, we show that every Hausdorff separable radial space is
R-separable and note that neither separable sequential nor separable
Whyburn spaces have to be selectively separable. A space is called
\emph{d-separable} if it has a dense $\sigma$-discrete subspace. We call
a space $X$ \emph{D-separable} if for every sequence of dense subspaces $(D_n :
n\in\omega)$ one can pick discrete subsets $F_n\subset D_n$ such
that $\bigcup_{n\in\omega}F_n$ is dense in $X$. Although $d$-separable spaces are often also $D$-separable (this is the case, for example, with linearly ordered $d$-separable or stratifiable spaces), we offer three examples of countable non-$D$-separable spaces. It is known
that d-separability is preserved by arbitrary products, and that for
every $X$, the power $X^{d(X)}$ is d-separable. We show that D-separability is
not preserved even by finite products, and that for every infinite
$X$, the power $X^{2^{d(X)}}$ is not D-separable. However, for every $X$ there
is a $Y$ such that $X\times Y$ is D-separable. Finally, we discuss
selective and D-separability in the presence of maximality. For example, we
show that (assuming ${\mathfrak d}=\mathfrak c$) there exists a maximal
regular countable selectively separable space, and that (in ZFC)
every maximal countable space is D-separable (while some of those
are not selectively separable). However, no maximal space satisfies the natural game-theoretic strengthening of D-separability.
\end{abstract}

\maketitle

\section{Introduction}
The area known as \emph{Selection principles in Mathematics} deals with selective variations of classical topological notions like compactness or separability (see \cite{Boazsurvey} or \cite{Marionsurvey} for a survey and \cite{Boazproblemsurvey} for another survey concentrating on open problems in the field). New results, questions and papers in the area are announced on the periodical SPM bulletin \cite{SPMbulletin}. Looking at the selective version of a certain property adds a combinatorial skeleton to it that often makes it easier to deal with. For example, Leandro Aurichi \cite{Leandro} has recently given one of the few known partial solutions to Eric Van Douwen's evasive \emph{D-space problem} (see \cite{Eisworth}) by replacing the Lindel\"of property with one of its selective strengthenings, the Menger property.

In this paper we will be concerned with the notion of selective separability and its variations. This notion has gained particular attention recently, as witnessed by the papers \cite{Marion1},
\cite{DiMaioKocinacMeccarielo}, \cite{BBMT}, \cite{Variations},
\cite{Addendum}, \cite{Babinkostova}, \cite{BermanDowTopProc},
\cite{RepovsZdomskyy}, \cite{GruenhageSakai}. A space $X$ is {\it
selectively separable} (also called {\it M-separable} or {\it SS})
if for every sequence of dense subspaces $(D_n : n\in\omega)$ one
can pick finite $F_n \subset D_n$ so that $\bigcup_{n\in\omega} F_n$
is dense in $X$. $X$ is {\it H-separable} if for every sequence of
dense subspaces $(D_n : n\in\omega)$ one can pick finite $F_n
\subset D_n$ so that every non-empty open set in $X$ intersects all
but finitely many $F_n$. $X$ is {\it R-separable} if for every
sequence of dense subspaces $(D_n : n\in\omega)$ one can pick  $p_n
\in D_n$ so that $\{p_n : n\in\omega\}$ is dense in $X$. $X$ is {\it
GN-separable} if $X$ is crowded and for every sequence of dense
subspaces $(D_n : n\in\omega)$ one can pick  $p_n \in D_n$ so that
$\{p_n : n\in\omega\}$ is {\it groupable}. This means that one can
find pairwise disjoint non-empty and finite sets $A_m$ for $m<\omega$
in such a way that $\{p_n : n\in\omega\} = \bigcup\{A_m : m\in\omega\}$
and every non-empty open set in $X$ intersects all but finitely many
$A_m$.

$X$ is SS$^+$ if Two has a winning strategy in the following game
${\sf G}_{\rm fin}({\mathcal D},{\mathcal D})$. $\mathcal D$ is the
collection of all dense subspaces of $X$. One picks $D_0\in \mathcal
D$, then Two picks a finite $F_0 \subset D_0$, then One picks $D_1
\in \mathcal D$, etc. Two wins if $\bigcup_{n\in\omega} F_n$ is
dense in $X$. (The term SS$^+$ is from \cite{BermanDowTopProc} but
the notion and the game ${\sf G}_{\rm fin}({\mathcal D},{\mathcal
D})$ were introduced in \cite{Marion1}.)
Barman and Dow discovered \cite{BermanDowTopProc} that every
separable Fr\'echet space is selectively separable. Gruenhage and Sakai \cite{GruenhageSakai} pointed out that separable Fr\'echet spaces are
even R-separable, and, if there are no isolated points,
GN-separable. In Section~\ref{ConvergenceSection} we discuss
the possibility to extend these results to spaces
satisfying convergence-type conditions weaker than Fr\'echet. It turns out that every regular separable radial space is
selectively separable while separable sequential spaces or separable
Whyburn spaces need not be selectively separable. We also consider the special case of countably compact spaces.

In Section~\ref{DSeparabilitySection} we consider a weaker form of
selective separability: the sets $F_n$ are supposed to be discrete
rather than finite; we call this property {\it D-separability}. This may be also viewed as a natural selective strenghtening of the notion of $d$-separability. Recall that $X$ is called {\it d-separable}
\cite{AVBicompactaAnd} (see also \cite{AVdSeparability},
\cite{Amirjanov}, \cite{Shap}, \cite{Tka2005}, \cite{JSz}) if $X$
has a dense $\sigma$-discrete subspace. The notion of $d$-separability is almost as old as separability and was introduced by Kurepa in his Ph.D. dissertation (see also \cite{Kurepa}), where it is called \emph{condition $K_0$}.

 It turns out that in some
cases d-separable spaces are D-separable. However, the behavior of
d-separability and D-separability, is quite different. This is particularly apparent if one looks at the product operation. Every product of d-separable spaces is
d-separable \cite{AVdSeparability}; for every T$_1$ space $X$, a
high enough power of $X$ is d-separable \cite{JSz} while we show
that there are two D-separable spaces with a non-D-separable
product, and every (Tychonoff) space has some power which is not
D-separable.

In Section~\ref{MiscellaniaSEction}, we discuss selective separability and
D-separability in maximal spaces. For example, we show that
(assuming ${\mathfrak d}=\mathfrak c$) there exists a maximal countable
selectively separable space, and that (in ZFC) every maximal regular
countable space is D-separable (while some of those are not
selectively separable). However, no maximal space is D$^+$-separable
(D$^+$-separability is a property stronger than D-separability and
defined in terms of topological games, see
Definition~\ref{DefinitionOfDPlusSep} below).

\section{Terminology and preliminaries}
For undefined topological notions we refer to \cite{Eng}, while for undefined set-theoretic notions we refer to \cite{Jech}. The letter $X$
always denotes a topological space. $X$ is {\em Fr\'echet} if for
every non-closed $A\subset X$ and every $p\in \overline{A}\setminus
A$ there is a sequence from $A$ converging to $p$. $X$ is {\it
sequential} if whenever $A$ is non-closed there are a $p\in
\overline{A}\setminus A$ and a sequence from $A$ converging to $p$.
$X$ has {\it countable tightness} if whenever $p\in \overline{A}$
there is a countable $B\subset A$ such that $x\in\overline{B}$. $X$
has {\it countable fan tightness} \cite{AVHurewiczFanTightness} if
whenever $p\in \overline{A_n}$ for all $n\in \omega$ one can pick
finite $F_n \subset A_n$ so that $p\in
\overline{\bigcup_{n\in\omega} F_n}$. $X$ has {\it countable strong
fan tightness} \cite{Sakai88} if whenever $p\in \overline{A_n}$ for
all $n\in \omega$ one can pick  $p_n \in A_n$ so that $p\in
\overline{\{p_n : n\in \omega\}}$. $X$ has {\it dense fan tightness}
if $X$ satisfies the definition of fan tightness restricted to $A_n$
dense in $X$. We will be using the following simple proposition without explicit mention.

\begin{proposition}\label{PropositionFromPreliminaries}
Let $X$ be separable. Then:

(1) {\rm \cite{BBMT}} $X$ is selectively separable iff for every
$p\in X$ and every sequence $(D_n : n\in\omega)$ of dense subspaces
of $X$ one can pick finite sets $F_n \subset D_n$ so that $p\in
\overline{\bigcup_{n\in\omega} F_n}$ (in other words, $X$ has
countable tightness with respect to dense sets).

(2) $X$ is H-separable iff for every $p\in X$ and every sequence
$(D_n : n\in\omega)$ of dense subspaces of $X$ one can pick finite sets $F_n \subset D_n$ so that every neighborhood of $p$ meets all but
finitely many $F_n$.

(3) $X$ is R-separable iff for every $p\in X$ and every sequence
$(D_n : n\in\omega)$ of dense subspaces of $X$ one can pick points $p_n
\in D_n$ so that $p\in \overline{\{p_n : n\in\omega\}}$ (in
other words, $X$ has countable strong tightness with respect to
dense sets).

(4) $X$ is GN-separable iff for every $p\in X$ and every sequence
$(D_n : n\in\omega)$ of dense subspaces of $X$ one can pick points $p_n
\in D_n$ and represent $\{p_n : n\in \omega\} =
\bigcup_{m\in\omega} A_m$ where the sets $A_m$ are non-empty,
finite, and pairwise disjoint, so that every neighborhood of $p$
intersects all but finitely many $A_m$.
\end{proposition}

$\delta(X)=\sup\{d(D): D$ is dense in $X$. If $X$ is compact, then
$\delta(X)=\pi w(X)$ \cite{JShelah}. Obviously, $\delta(X)=\omega$
for every selectively separable space$X$.

$X$ is {\it radial} if for every $A\subset X$ and every
$p\in\overline{A}$ there is a well-ordered net $\{x_\alpha :
\alpha<\kappa\}\subset A$ which converges to $p$. $X$ is {\it
pseudoradial} if for every non-closed $A\subset X$ there is a
$p\in\overline{A}\setminus A$ and a well-ordered net $\{x_\alpha :
\alpha<\kappa\}\subset A$ which converges to $p$. A set $A\subset X$
is $\kappa$-{\it closed} (where $\kappa$ is a cardinal) if
$\overline{B}\subset A$ whenever $B\subset A$ and $|B|\leq\kappa$.
$X$ is {\it semiradial} (see \cite{BellaMoreOnSequentialProperties},
\cite{BellaTironiEncyclopedia}) if for every $\kappa$, every
non-$\kappa$-closed set $A$ contains a well-ordered net of length
$\leq\kappa$ converging to a point outside $A$. Among the various
subclasses of pseudoradial spaces considered in the literature, the
class of semiradial spaces is the smallest one which includes all
radial and all sequential spaces.

$X$ has the {\it Whyburn property} if for every $A\subset X$ and
every $p\in\overline{A}\setminus A$ there exists $B\subset A$ such
that $\overline{B} = A \cup \{p\}$. Every Fr\'echet space is Whyburn and every compact Whyburn space is Fr\'echet (see \cite{TkaYas}). For $p \in \omega^*$, the space $\omega \cup \{p\}$ with the topology inherited from $\beta \omega$ is a non-Fr\'echet Whyburn topological space. The space $C_p([0,1])$ of all
continuous functions from $[0,1]$ to $\mathbb{R}$ with the topology
of pointwise convergence is a nice Whyburn topological group which is not
Fr\'echet \cite{BellaYaschenkoAPandWAP}.

Recall that for $M\subset X$, ${\rm seqcl}(M) = \{x\in X :$ there is
a sequence converging from $M$ to $x\}$ and ${\rm seqcl}_\alpha(M)$
is defined inductively by ${\rm seqcl}_\alpha(M)= {\rm
seqcl}(\bigcup_{\beta<\alpha} {\rm seqcl}_\beta(M))$. If $X$ is
sequential then there exists an ordinal $\alpha^*$ called the {\it
sequential order} of $X$ such that ${\rm seqcl}_{\alpha^*}(M) =
\overline{M}$ for every $M\subset X$. The sequential order of any
sequential space is $\leq\omega_1$.

Let ${}^n\omega$ be the set of all  functions $s:n=\{0,1, \ldots,
n-1\}\to \omega$ and let $Seq=\bigcup\{{}^n\omega: n<\omega\}$. If
$s\in {}^n\omega $ and $k\in \omega$,  we write $s\frown
k=s\cup\{(n,k)\}\in {}^{n+1}\omega$.  Given a free filter $\mathcal
F$ on $\omega$, we denote by   $Seq(\mathcal F)$  the topological
space having $Seq$ as the underlying set and  the topology obtained by
declaring a set $U\subseteq Seq$ open if and only if for any $s\in
U$  $\{n: s\frown n \in U\}\in \mathcal F$. $Seq(\mathcal F)$ is
always a Hausdorff zero-dimensional dense-in-itself space (see
\cite{VaughanSeqP} for more information). In several instances in
this paper, we will use $Seq({\mathcal F})$ where $\mathcal F =
\{A\subset \omega : |\omega\setminus A|<\omega\}$ is the Fr\'echet
filter. This space is also known under the name $S_\omega$
\cite{ArhFranklin}. $Seq({\mathcal F})$ is sequential of sequential
order $\omega_1$.

$X$ is {\it crowded} (also called {\it dense in itself}) if $X$ does
not have isolated points. $X$ is {\it maximal} if $X$ is crowded and
no topology strictly stronger than the topology of $X$ is crowded.
$X$ is {\it resolvable} ($\omega$-{\it resolvable}) if $X$ contains
two (respectively, a countably infinite family of pairwise) disjoint
dense subspaces. $X$ is {\it submaximal} if every subset is open in
its closure, or, equivalently (see \cite{ArhCollinsSubmaximal}) if
the complement of every dense set is closed and discrete. Every
maximal space is submaximal. Every crowded submaximal (hence every
maximal) space is irresolvable (= not resolvable). $X$ is {\it
Baire} if no non-empty open set in $X$ is representable as the union
of countably many nowhere dense sets. $X$ is {\it strongly
irresolvable} \cite{KunenSzymanskiTall} if all non-empty open sets
are irresolvable. Strongly irresolvable Baire is abbreviated as SIB
\cite{KunenSzymanskiTall}. For any space $X$, the \emph{dispersion character} $\Delta(X)$ of $X$ is defined as the
minimum cardinality of a non-empty open set in $X$. A crowded space
$X$ is {\it extra-resolvable} if there is a family $\mathcal G$ of
dense subspaces of $X$ such that $|{\mathcal G}|> \Delta(X)$ and for
every two distinct $G, G^\prime \in \mathcal G$, $G\cap G^\prime$ is
nowhere dense.

$X$ is {\it discretely generated} \cite{DTTW} if whenever $p\in
\overline{A}$ there is a discrete $D\subset A$ such that
$p\in\overline{D}$.

$X$ is a $\sigma$-{\it space} if $X$ has a $\sigma$-discrete
network.

$X$ is {\it monotonically normal} if one can assign to
every point $x \in X$ and open set $U \subset X$ an open set $H(x,U) \subset U$ such that $x \in H(x,U)$ and if $H(x,U) \cap H(y,V) \neq \emptyset$ then either $x \in V$ or $y \in U$. The function $H$ is called a \emph{monotone normality operator}.

$X$ is {\it stratifiable} if one can assign to every $n\in\omega$ and every
closed set $H\subset X$ an open set $G(n,H)\supset H$ so that $H=
\bigcap_{n\in\omega} \overline{G(n,H)}$ and $G(n,H) \subset G(n,K)$
whenever $H\subset K$. Every stratifiable space is both
monotonically normal and a $\sigma$-space \cite{GaryHandbook}.

A set $D \subset X^2$ is called {\it slim}
\cite{GruenhageNatkaniecPiotrowskiThinVeryThinSlim} if the
intersection of $D$ with every {\it cross-section} $(\{p\} \times X)
\cup (X \times \{p\})$ is nowhere dense (in this cross section). If
$B\subset A$ then $\pi_B : \prod\limits_{\alpha\in A}X_\alpha \to
\prod\limits_{\alpha\in B}X_\alpha$ denotes the projection of the
product onto a subproduct.

$cov({\mathcal M})$ is the minimum cardinality of a family of
nowhere dense subsets of $\mathbb R$ that covers $\mathbb R$. A function
$g\in \omega^\omega$ is said to {\it guess} the family of functions $\Phi
\subset \omega^\omega$ if for every $f\in \Phi$, $f(n)=g(n)$ for
infinitely many $n$. It is known that if $|\Phi| < cov({\mathcal
M})$, then there is $g$ that guesses $\Phi$, see
\cite{BartosinskyJudah}.

\section{Convergence and selective
separability}\label{ConvergenceSection}

Gruenhage and Sakai \cite{GruenhageSakai} observed that separable
Fr\'echet spaces are R-separable. Basically, there are three
\lq\lq natural ways" to try to strengthen this result: one is to
move from Fr\'echet to radial, another is to move from
Fr\'echet to sequential (or, more generally, to spaces of
countable tightness), and yet another is from Fr\'echet to
Whyburn.

\subsection{Radial spaces}

\begin{proposition}
(1) A Hausdorff separable radial (with respect to dense subspaces)
space $X$ is R-separable.

(2) If, in addition, $X$ does not have isolated points, then $X$ is
GN-separable.
\end{proposition}

\begin{proof}
(1) It suffices to show that, given a point $p\in X$
and a sequence of dense subspaces $(D_n : n\in\omega)$ one can pick
$p_n \in D_n$ so that $p\in \overline{\{p_n : n\in\omega\}}$. Assume
$p$ is not isolated, otherwise $p$ is contained in every dense set and the statement we want to prove becomes trivial. Let $\mathcal U$ be a
maximal pairwise disjoint family of non-empty open sets in $X$ such
that $x\not\in\overline{U}$ for every $U\in \mathcal U$ (*). Since
$X$ is Hausdorff and $p$ is not isolated, $\bigcup\mathcal U$ is
dense in $X$. Since $X$ is separable, $\mathcal U$ is countable;
enumerate it as $(U_n : n\in\omega)$. Put $Y= \bigcup\{D_n \cap U_n
: n\in\omega\}$. Then $Y$ is dense in $X$ and thus there is
$S\subset Y$ which can be enumerated as a well-ordered net
converging to $p$. We can assume that $S$ is of minimal cardinality among all well-ordered nets contained in $Y$ which converge to $p$,
so $|S|$ is regular, and then it follows from (*) and the
countability of $\mathcal U$ that $S$ is a convergent sequence.
Again from (*), $S$ must have non-empty intersection with infinitely
many sets $D_n\cap U_n$; pick a subsequence $S^\prime \subset S$
that intersects each $D_n\cap U_n$ at at most one point.

For $n\in\omega$, if $S^\prime \cap D_n\cap U_n$ is non-empty, let
$p_n$ be the unique point in this intersection. Otherwise pick $p_n$
arbitrarily. Then the points $p_n$ are as desired.

(2) Partition $\{p_n : n\in\omega\}$ from part 1 into pairwise
disjoint finite sets $A_m$ so that each $A_m$ contains at least one
point of $S^\prime$ and apply
Proposition~\ref{PropositionFromPreliminaries}, part 4.
\end{proof}

\vspace{3mm} The above proposition cannot be extended to pseudoradial spaces
because as we will see below even separable sequential spaces need not be selectively separable.

\begin{corollary} \label{corradial}
Every compact separable radial space has countable $\pi$-weight.
\end{corollary}

\begin{proof}
For a selectively separable $X$,
$\delta(X)=\omega$, and for a compact $X$, $\delta(X)= \pi w(X)$
\cite{JShelah}.
\end{proof}

\begin{corollary} \label{corMN}
Every separable compact monotonically normal space has countable $\pi$-weight.
\end{corollary}

\begin{proof}
Monotonically normal spaces are radial \cite{WilliamsZhouMonotNormal}.
\end{proof}

\vspace{3mm} While we could not find a reference for
Corollary~\ref{corradial}, Corollary~\ref{corMN} can be also derived
from \cite{Gartside}, Corollary~19 (which says that density equals
$\pi$-weight for monotonically normal compact spaces).

One can wonder whether Corollary $\ref{corradial}$ can be extended to semiradial spaces. Some mild evidence is provided by the fact that every compact sequential separable space has countable $\pi$-weight (this follows easily from the inequality $\pi \chi(X) \leq t(X)$, for every compact space $X$ see \cite{J}) and every sequential space is semiradial. However, the answer is consistently negative. Bella \cite{BellaMoreOnSequentialProperties} showed
that the space $2^{\omega_1}$ is semiradial if and only if ${\mathfrak
p} > \omega_1$. One might hope at least for a consistency result. Indeed, A. Dow proved
\cite{DowCompactSeparableRadial} that there are models in which
every compact separable radial space is Fr\'echet, but it is still
an open problem if there are models in which every compact separable
semiradial space is sequential. In view of
Corollary~\ref{corradial}, we suggest a weaker form of this problem.

\begin{question}
{\rm Is it consistent that every compact separable semiradial space
has countable $\pi$-weight?
 }
\end{question}

\subsection{Sequential spaces and spaces of countable tightness}


Strengthening the Barman-Dow result in this direction is in
general not possible.

If $K$ is an infinite subset of $\omega$, then it is easy to check
that the set $D_K=\bigcup\{{}^k\omega:k\in K\}$ is dense in
$Seq(\mathcal F)$. Moreover, it is also quite easy to realize that,
for any choice of a finite set $F_n \subseteq {}^n\omega $, the set
$\bigcup\{F_n:n<\omega\}$ is closed and nowhere dense in
$Seq(\mathcal F)$.   Taking this into account, we see that if
$\{H_n=\bigcup\{{}^k\omega:n\le k<\omega\}$, then  the sequence of
dense sets  $\{H_n:n<\omega\}$ witnesses that the space
$Seq(\mathcal F)$ is never selectively separable.

If as $\mathcal F$ we take the filter of all cofinite subsets of
$\omega$, then  $Seq(\mathcal F)$  turns out to be sequential.
Indeed, if $A$ is a non-closed subset of $Seq(\mathcal F)$, then
$Seq(\mathcal F)\setminus A$ is not open and so there is some $s\in
Seq(\mathcal F)\setminus A$ such that    the set $\{n:s\frown n\in
Seq(\mathcal F)\setminus A\}$  has an infinite complement $E$.
Therefore, $S=\{s\frown n:n\in E\}\subseteq A$ and we immediately
 check that $S$ converges to $s$.

Thus, there is a countable sequential space which is not selectively
separable.   However, this space has sequential order $\omega_1$.
So, what is left open is:

\begin{question}
{\rm Is every Hausdorff separable sequential space of finite or
countable sequential order selectively separable?}
\end{question}

That separable spaces of countable tightness need not be
selectively separable is well known. There are many examples such as
$C_p($Irrationals$)$ or even some countable spaces \cite{BBMT},
\cite{Variations}. However, adding some restrictions on the character of points or some covering properties we can get positive results. Here are a few of the most interesting.

\begin{proposition} \ \label{ThreeCasesCountableTightness}

\begin{enumerate}
\item {\rm \cite{BBMT}} If a separable space $X$ of countable
tightness has a dense set of points of character less than $\mathfrak
d$, then $X$ is selectively separable.

\item \label{itemcountablycompact} A regular countably compact separable space of countable
tightness is selectively separable.

\item \label{itemGruenhageSakai} {\rm \cite{GruenhageSakai}} More generally, let $X$ be a regular
separable space of countable tightness. If each point is contained
in a countably compact set of countable character in $X$, then $X$
is selectively separable.
\end{enumerate}
\end{proposition}

Item $\ref{itemcountablycompact}$ follows directly from the fact that a regular countably
compact space of countable tightness has countable fan tightness
\cite{ArhBellaCountableFanTightness}. In order to slightly improve
this result we prove a proposition which may be of some independent interest.

\begin{proposition}\label{CountablyCompactStrongFanTightness}
A regular countably compact space $X$ of countable tightness has
countable strong fan tightness.
\end{proposition}

\begin{proof}
Let $b\in X$, $A_n\subset X$, $b\in \overline{A_n}$ for
all $n\in\omega$. Without loss of generality we assume that the sets
$A_n$ are countable. Put $\tilde{X}= \overline{\bigcup_{n\in\omega}
A_n}$. Then $\tilde{X}$ is a regular separable countably compact
space of countable tightness; $b\in \tilde{X}$. Fix a base $\mathcal
B$ of neighborhoods of $b$ in $\tilde{X}$ such that $|{\mathcal
B}|\leq \mathfrak c$. For every $U\in \mathcal B$ fix an open in
$\tilde{X}$ set $V(U)$ such that $b\in V(U) \subset \overline{V(U)}
\subset U$. Last, fix an almost disjoint family $\mathcal R$ of
infinite subsets of $\omega$ enumerated by $\mathcal B$: ${\mathcal
R} = \{N_U : U\in{\mathcal B}\}$.

Let $U\in \mathcal B$. For every $n\in N_U$ pick $x_{n,U} \in A_n
\cap V(U)$. Since $\tilde{X}$ is countably compact the set
$\{x_{n,U} : n\in N_U\}$ has a limit point, say $s_U$. Then $s_U \in
\overline{V(U)} \subset U$.

Put $S=\{s_U : U\in {\mathcal B}\}$. Then $b\in\overline{S}$. Since
$\tilde{X}$ has countable tightness there is a countable subfamily
${\mathcal B}_0 \subset \mathcal B$ such that $b\in\overline{S_0}$
where $S_0 = \{s_U : U\in {\mathcal B}_0\}$. Enumerate ${\mathcal
B}_0 = \{U_m : m\in \omega\}$. For $U_m \in {\mathcal B}_0$ put
$\tilde{N}_{U_m} = N_{U_m} \setminus \bigcup_{k<m} N_{U_k}$. Then
the sets $\tilde{N}_{U_m}$ are pairwise disjoint, and
$\tilde{N}_{U_m}$ differs from $N_{U_m}$ only in finitely many
points.

Let $n\in \omega$. If $n\in \tilde{N}_{U_m}$ for some (single!) $m$
then put $a_n = x_{n,U_m}$. Otherwise pick $a_n \in A_n$
arbitrarily. Thus we have $a_n \in A_n$ defined for all $n$, and
$b\in \overline{\{a_n : n\in\omega\}}$ (because $b\in\overline{S_0}$
and each point of $S_0$ is in $\overline{\{a_n : n\in\omega\}}$).
\end{proof}

\vspace{3mm} As an immediate corollary we get that regular countably
compact separable spaces of countable tightness are R-separable.
However, we are going to prove a stronger result which is a
simultaneous improvement of all parts of
Proposition~\ref{ThreeCasesCountableTightness}. Let $X$ be a space
and $x\in X$. Weakening the definition of the cardinal function
$h(x,X)$ \cite{Eng} we denote by $h^*(x,X)$ the smallest cardinal
number $\kappa$ such that there exists a countably compact $H\subset
X$ with $x\in H$ and $\chi(H,X)=\kappa$.

\begin{theorem} \label{generalizationofGruenhageSakai}
Let $X$ be a regular  separable space of countable tightness.

(a) If the inequality $h^*(x,X)<\mathfrak d$ holds in a dense set  of
points, then $X$ is selectively  separable.

(b) If the inequality $h^*(x,X)< cov(\mathcal M)$ holds in a dense
set of points, then $X$ is R-separable.
\end{theorem}

\begin{proof}
Without any loss of generality, we may suppose that $X$
does not have isolated points. Let $(D_n:n\in \omega)$ be a sequence
of dense subsets of $X$. Since from the hypotheses it follows that
any dense set is separable,  we may assume that each $D_n$ is
countable.

\underline{Part a.} Let us begin by fixing a countable dense set $C$
such that $h^*(x,X)<\mathfrak d$ holds for each $x\in C$. Let
$\{J_x:x\in C\}$ be a partition of $\omega$ in infinite sets and fix
for each $x\in C$ a countably compact subspace $H_x$ satisfying
$x\in H_x$ and $\chi(H_x,X)=\kappa_x <\mathfrak d$. Let $\Phi=\prod
_{n\in J_x} [D_n]^{<\omega}$ and for each $\phi\in \Phi$ let
$[\phi]= D({\bigcup\{\phi(n):n\in J_x\}})$ (here $D(S)$ indicates
the derived set of $S$).

\medskip

\noindent {\bf Claim}. $x\in \overline {\bigcup\{[\phi]:\phi\in \Phi\}}$.

\begin{proof}[Proof of Claim]  Fix an open neighborhood $V$ of $x$. We are
going to check that $\overline V\cap \bigcup\{[\phi]:\phi\in
\Phi\}\ne \varnothing$. \underline{Case 1} If there exists an
infinite set $J\subseteq J_x$ such that $V\cap D_n\cap H_x $ is
infinite for each $n\in J$, then we may define a function $\psi\in
\Phi$  by letting $\psi(n)=\{p_n\}$, for some  $p_n\in V\cap D_n\cap
H_x$ if $n\in J$, and $\psi(n)=\varnothing$ otherwise.  The function
$\psi$ can be made one-to-one on $J$ and so $\bigcup\{\psi(n):n\in
J_x\}$ is an infinite subset of the countably compact set $H_x$.
This guarantees that $[\psi]\ne \varnothing $ and clearly by
construction $[\psi]\subseteq \overline V$. \underline{Case 2} If
such a $J$ does not exist, then we may  assume, without any loss of
generality, that $V\cap D_n\cap H_x=\varnothing$ for each $n\in J_x$
(it suffices to replace $D_n$ with $D_n \setminus (V\cap H_x)$ for
all but finitely many $n\in J_x$ and then remove a finite part of
$J_x$). For any $n\in J_x$ write $D_n =\{x_{n,k}: k< \omega\}$ and
let $\{U_\alpha :\alpha \in \kappa_x\}$ be a local base of $H_x$ in
$X$. For any $\alpha$ fix a function $f_\alpha :J_x\to \omega$,
defined in such a way that $x_{n,f_\alpha (n)}\in V\cap U_\alpha$
for each $n\in J_x$. Since $\kappa_x <\mathfrak d$, there is a function
$g:J_x\to \omega$ such that for each $\alpha $ the inequality
$f_\alpha (n)\le g(n)$ holds for infinitely many $n$'s. Let $\psi$
be the element of $\Phi$ defined by  letting $\psi(n)=
\{x_{n,k}:k\le g(n)\}$. By the choice of $g$, for each $\alpha $
there is some $n\in J_x$ such that $U_\alpha \cap \psi(n)\ne
\varnothing$ and this in turn implies that we must have
$\overline{\bigcup \{\psi(n):n\in J_x\}}\cap H_x\ne \varnothing$.
Since the values of $\psi$ are disjoint from $H_x$, the latter
formula implies that $[\psi]\ne \varnothing$. Of course, by
construction we have $[\psi]\subseteq
 \overline V$ and the claim is proved.
 \renewcommand{\qedsymbol}{$\triangle$}
 \end{proof}

Now, thanks to the Claim and the countable tightness of $X$, there
are countably many $\phi_n\in \Phi$ and points $z_n\in [\phi_n]$
such that $x\in \overline {\{z_n:n\in \omega\}}$.  Observe that if
$O$ is an open neighborhood of $x$ then $z_n\in O$ for some $n$ and,
being $z_n$ an accumulation point of the set $\bigcup\{\phi_n(k):
k\in J_x\}$, we actually have $O\cap \bigcup\{\phi_n(k):k\in J_x,
k\ge n\}\ne\varnothing$. Now, we put $F^x_n=\bigcup\{\phi_n(k):k\le
n\}$. As in the first paragraph, we have $x\in \overline {\bigcup
\{F^x_n:n\in J_x\}}$.  To finish, it is enough to note that  the
chosen family $\{F^x_n: n\in J_x, x\in C\}$ is dense in $X$.

\underline{Part b.} Fix a countable dense set $C$ such that
$h^*(x,X)<cov({\mathcal M})$ holds for each $x\in C$.  Let
$\{J_x:x\in C\}$ be a partition of $\omega$ in infinite sets and fix
for each $x\in C$ a countably compact subspace $H_x$ satisfying
$x\in H_x$ and $\chi(H_x,X)=\kappa_x <cov({\mathcal M})$.

Let $x\in C$. Fix a base $\mathcal B$ of neighborhoods of $x$ in $X$
such that $|{\mathcal B}|\leq\mathfrak c$. Assign to every $U\in
\mathcal B$ an open neighborhood $V_U$ of $x$ so that $x\in V_U
\subset \overline{V_U} \subset U$. Fix an almost disjoint family
$\{N_{x,U} : U\in{\mathcal B}\}$ of infinite subsets of $\omega$.
\medskip

\noindent {\bf Claim.} For every $U\in \mathcal B$ there is $\phi \in
\prod_{n\in N_{x,U}} D_n$ such that $\overline{V_U} \cap D(\{\phi(n)
: n\in N_{x,U}\}) \neq\varnothing$.

\begin{proof}[Proof of Claim] In the case when there exists an infinite
$J\subset N_{x,U}$ such that $V_U \cap D_n \cap H_x$ is infinite for
each $n\in J$ the argument repeats the similar one from part a. So
assume no such $J$ exists and then without loss of generality we
assume that $V_U \cap D_n \cap H_x = \varnothing$ for for all $n\in
N_{x,U}$. Enumerate $D_n=\{x_{n,k} : k\in\omega\}$. Let $\{O_\alpha
: \alpha<\kappa_x\}$ (where $\kappa_x<cov({\mathcal M})$ be a local
base of $H_x$ in $X$. For every $\alpha$ fix a function $f_\alpha:
N_{x,U} \to \omega$ defined in such a way that $x_{n,f_\alpha(n)}
\in V_U \cap O_\alpha$ for each $n\in N_{x,U}$. Since $\kappa_x <
cov({\mathcal M})$ there is a function $g: N_{x,U} \to \omega$ such
that for each $\alpha$ the equation $f_\alpha(n)=g(n)$ holds for
infinitely many $n$'s. Define $\phi$ by $\phi(n) = x_{n,g(n)}$ for
each $n\in N_{x,U}$. This proves the claim.
\renewcommand{\qedsymbol}{$\triangle$}
\end{proof}

Now, using claim, for every $U\in \mathcal B$ fix $\phi_U \in
\prod_{n\in N_{x,U}} D_n$ and $z_U \in \overline{V_U} \cap
D(\{\phi(n) : n\in N_{x,U}\})$. Put $Z=\{z_U : U\in {\mathcal B}\}$.
Then $x\in \overline{Z}$ and, since $X$ has countable tightness,
there is a countable subset, say $Z_0 = \{z_{U_k} : k\in\omega\}$,
such that $x\in \overline{Z_0}$. For $k\in\omega$, put
$\tilde{N}_{x,U_k} = N_{x,U_k} \setminus \bigcup_{l<k} N_{x,U_l}$.
Then the sets $\tilde{N}_{x,U_k}$ ($k\in\omega$) are pairwise
disjoint and differ from $N_{x,U_k}$ only by finitely many elements.
Let $n\in J_x$. If $n$ belongs to some (then only to one)
$N_{x,U_k}$ then put $a_n = x_{n,g(n)}$. If not, choose $a_n\in D_n$
arbitrarily. Then $x\in\overline{\{a_n : n\in J_x\}}$.
\end{proof}

\vspace{3mm} A crucial role in the proof of Theorem
$\ref{generalizationofGruenhageSakai}$ as well as in Gruenhage and
Sakai's proof of Proposition $\ref{ThreeCasesCountableTightness}$,
part $\ref{itemGruenhageSakai}$ \cite{GruenhageSakai}  is played by
Proposition~\ref{ThreeCasesCountableTightness}, part
$\ref{itemcountablycompact}$ (and its variation,
Proposition~\ref{CountablyCompactStrongFanTightness}). This suggests
to look for some possible generalization. In one direction we may
try to weaken \lq\lq regular" to Hausdorff and in the other to
weaken countably compact to pseudocompact. Unfortunately, we have an
answer only for the first case.

 It is well known (see
\cite{VaughanSeqP}) that if $\mathcal U$ is a free ultrafilter on
$\omega$ then the space $Seq(\mathcal U)$ is extremally
disconnected. So, $X=Seq(\mathcal U)$ is a countable Hausdorff
zero-dimensional extremally disconnected non-selectively separable
space. Now, consider the \v Cech-Stone compactification $\beta X$ of
this $X$. Theorem~1.1g  of \cite{VaughanCountablyCompactT2} shows
that there exists a strengthening of the topology of $\beta X$ in
such a way that the resulting space $Y$ has the following
properties:
\begin{enumerate}
\item  $X$ is a dense subspace of $Y$;
\item  $Y$ is locally countable;
\item  each closed infinite subset of $Y$  has cardinality $2^{\mathfrak c}$.
\end{enumerate}
So, we get:

\begin{example}
{\rm There exists a separable countably compact  Urysohn space  of
countable tightness which is not selectively separable.
}
\end{example}

Moving from countably compact to pseudocompact appears much harder.
Indeed,  with a lot of effort, Bella and Pavlov
\cite{BellaPavlovEmbeddingsPseudocompact} constructed a
Tychonoff pseudocompact space of countable tightness which does not
have countable fan tightness. But such a space has a countable set
of isolated points and so it is selectively separable. For these
reasons, the next problem sounds very interesting:

\begin{problem}
{\rm Find a Tychonoff pseudocompact separable space of countable
tightness which is not selectively separable.
}
\end{problem}

Forgetting separability, we may formulate a possibly easier problem:

\begin{problem} \label{pseudocompact}
{\rm Find a Tychonoff pseudocompact space of countable tightness
which does not have countable dense fan tightness.
}
\end{problem}

\subsection{Whyburn spaces}

Barman and Dow constructed a countable regular maximal space which
is not selectively separable \cite{BermanDowTopProc}. On the other
hand, the first author and I. Yaschenko showed in
\cite{BellaYaschenkoAPandWAP} that every regular maximal space has
the Whyburn property. Therefore we get:

\begin{corollary}
There exists a countable regular Whyburn space which is not
selectively separable.
\end{corollary}


Tkachuk and Yaschenko \cite{TkaYas} proved that every countably compact Whyburn space is Fr\'echet. So countably compact Whyburn separable spaces are selectively separable. However pseudocompact Whyburn spaces need not be Fr\'echet \cite{PTTW}, not even if they have countable tightness \cite{BellaSimon}. So, also in view of Problem $\ref{pseudocompact}$, we have the following question.

\begin{question}
Suppose $X$ is a pseudocompact Whyburn separable space. Is $X$ selectively separable? What if $X$ has countable tightness?
\end{question}

\section{D-separability}\label{DSeparabilitySection}

A space is called $d$-separable if it contains a $\sigma$-discrete dense subspace. We introduce some selective version of
this property.

\begin{definition}\label{DefinitionOfDsep}
{\rm $X$ is {\it D-separable} if for every sequence of
dense subspaces $(D_n : n\in\omega)$ one can pick discrete sets $F_n
\subset D_n$ so that $\bigcup_{n\in\omega} F_n$ is dense in $X$.

\item $X$ is {\it DH-separable} if for every sequence of dense
subspaces $(D_n : n\in\omega)$ one can pick discrete $F_n \subset
D_n$ so that every non-empty open set in $X$ intersects all but
finitely many $F_n$.
 }
\end{definition}

Consider the following games on a space $X$ (as above, $\mathcal D$
denotes the collection of all dense subspaces of $X$). In the game
${\sf G}_{\rm dis}({\mathcal D},{\mathcal D})$, One picks $D_0\in
\mathcal D$, then Two picks a discrete $F_0 \subset D_0$, then One
picks $D_1 \in \mathcal D$, etc. Two wins if $\bigcup_{n\in\omega}
F_n$ is dense in $X$. The game ${\sf G}_{\rm dis,H}({\mathcal
D},{\mathcal D})$ is similar, only Two wins if every non-empty open
set in $X$ intersects all but finitely many $F_n$.

\begin{definition}\label{DefinitionOfDPlusSep}
{\rm $X$ is {\it D$^+$-separable} if Two has a winning
strategy in ${\sf G}_{\rm dis}({\mathcal D},{\mathcal D})$. Say that
$X$ is {\it DH$^+$-separable} if Two has a winning strategy in ${\sf
G}_{\rm dis,H}({\mathcal D},{\mathcal D})$.
 }
\end{definition}

The following implications (where, for example, SS denotes
selectively separable and D denotes D-separable) are
straightforward.

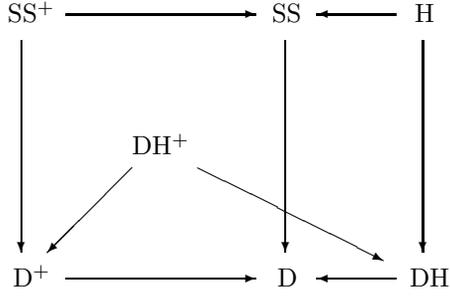
\begin{picture}(300,135)
\put(50,10){D$^+$} \put(150,10){D} \put(200,10){DH}
\put(70,13){\vector(1,0){72}} \put(195,13){\vector(-1,0){30}}
\put(95,60){DH$^+$} \put(95,55){\vector(-1,-1){32}}
\put(120,55){\vector(2,-1){70}}

\put(48,110){SS$^+$} \put(148,110){SS} \put(202,110){H}
\put(70,113){\vector(1,0){72}} \put(195,113){\vector(-1,0){30}}

\put(53,103){\vector(0,-1){80}} \put(153,103){\vector(0,-1){80}}
\put(205,103){\vector(0,-1){80}}
\end{picture}

We know that most of the arrows in the diagram cannot be reversed. To see that an arrow pointing from a selective separability-type property (top row of the diagram) to a selective $d$-separability-type property (bottow row and center of the diagram) cannot be reversed simply take any metric non-separable space. In some cases we will be able to improve this and obtain a separable counterexample. For instance, Example $\ref{answerboaz}$ is a countable space showing that $D \nrightarrow SS$. Any countable submaximal space shows that $D \nrightarrow D^+$ (see Corollary $\ref{ImmediateFromAngelosFive}$ and Theorem $\ref{D+separableomegaresolvable}$). However, the relationship between $D$ and $DH$ is not well-understood, and hence we have the following open problem

\begin{problem} \ {\rm
\begin{enumerate}
\item Find an example of a (countable) $D$-separable non-$DH$-separable space.
\item Find an example of a (countable) $D^+$-separable non-$DH^+$-separable space.
\item Find an example of a (countable) $DH$-separable, non-$DH^+$-separable space.
\end{enumerate}}
\end{problem}

\subsection{Which spaces are D-separable?}

\begin{proposition}\label{TrivialCasesOfDSeparability}
(1) Every space with a $\sigma$-disjoint $\pi$-base is
DH$^+$-separable.

(2) Every space with a $\sigma$-locally finite $\pi$-base is
DH$^+$-separable.

(3) Every T$_1$ space with a $\sigma$-closure preserving $\pi$-base
is DH$^+$-separable.
\end{proposition}
\begin{lemma}\label{ReductionDisjointFamilies}
If $\mathcal U$ and $\mathcal V$ are pairwise disjoint families of
non-empty sets in $X$, then there is a pairwise disjoint family
$\mathcal W$ of non-empty sets in $X$ such that:
\begin{enumerate}
\item Every element of $\mathcal U$ contains an element of $\mathcal
W$;
\item Every element of $\mathcal V$ contains an element of $\mathcal
W$;
\item Every element of $\mathcal W$ is contained in some element of
${\mathcal U} \cup {\mathcal V}$.
\end{enumerate}
\end{lemma}

\begin{proof}
Put ${\mathcal W}_0 = \{U\cap V : U\in{\mathcal U},
V\in{\mathcal V}$ and $U\cap V\neq\varnothing\}$ and ${\mathcal W} =
{\mathcal W}_0 \cup \{U\in {\mathcal U} :$ there is no $W\in
{\mathcal W}$ with $U\supset W\} \cup \{V\in {\mathcal V} :$ there
is no $W\in {\mathcal W}$ with $V\supset W\}$.
\end{proof}

\begin{lemma}\label{AngelosLemmaOnLocallyFinite}
If $\mathcal U$ is a locally finite family of non-empty open sets in
a space $X$ and $D$ is a dense subspace of $X$, then there is a
discrete set $A$ such that $A\cap U\neq \varnothing$ for every non-empty $U\in \mathcal U$.
\end{lemma}

\begin{proof}
\footnote{In the T$_1$ case, the proof is trivial: just
pick a point in every non-empty element of $\mathcal U$. But the
statement is valid without any assumption on separation.} For every
$U\in\mathcal U$ pick a point $p_U \in U\cap D$. Let $B=\{p_U :
U\in{\mathcal U}\}$. The local finiteness of $\mathcal U$ implies
that every $x\in B$ is contained in a set $V_x\subset B$ which is
finite and open in $B$. We choose as $V_x$ an open set of minimum
size. Let us call a point $x\in B$ good if for each $y\in V_x$ we
have $V_y = V_x$. It is clear that for each $x\in B$ there is a good
point $y$ such that $V_y \subset V_x$. Moreover, if $y$ and $z$ are
good, then either $V_y = V_z$ or $V_y \cap V_z = \varnothing$. now,
fix a well ordering on $B$ and for each good point $y\in B$ let
$a(y)= \min V_y$. The set $A$ of all such $a(y)$ is discrete and $A$
intersects every element of $\mathcal U$.
\end{proof}

\begin{proof}[Proof of Proposition~\ref{TrivialCasesOfDSeparability}]
(1) Let ${\mathcal U} = \bigcup_{n\in\omega}{\mathcal U}_n$ (where each
${\mathcal U}_n$ is pairwise disjoint and consists of non-empty
sets) be a $\pi$-base of $X$. Applying
Lemma~\ref{ReductionDisjointFamilies} inductively one gets pairwise
disjoint families ${\mathcal W}_n$ of non-empty open sets such that
whenever $m\leq n$, every element of ${\mathcal W}_n$ is contained
in some element of ${\mathcal U}_m$, and every element of ${\mathcal
U}_m$ contains an element of ${\mathcal W}_n$. At the $n$th inning
One chooses a dense subspace $D_n$ and Two picks for every $U\in
{\mathcal U}_n$ a point $p_{D_n,U} \in U$ and sets $F_n =
\{p_{D_n,U} : U\in{\mathcal U}_n\}$.

(2) Let ${\mathcal U} = \bigcup_{n\in\omega}{\mathcal U}_n$ (where
each ${\mathcal U}_n$ is locally finite) be a $\pi$-base of $X$. Put
${\mathcal W}_n = \bigcup_{m\leq n}{\mathcal U}_m$. At the $n$th
inning One chooses a dense subspace $D_n$ and Two uses
Lemma~\ref{AngelosLemmaOnLocallyFinite} to find a discrete subspace
$F_n \subset D_n$ which meets every element of ${\mathcal W}_n$.

(3) Let ${\mathcal U} = \bigcup_{n\in\omega}{\mathcal U}_n$ (where
each ${\mathcal U}_n$ is closure preserving) be a $\pi$-base of $X$.
At the $n$th inning One chooses a dense subspace $D_n$ and Two picks
for every $U\in {\mathcal U}_n$ a point $p_{D_n,U} \in U$ and sets
$F_n = \{p_{D_n,U} : U\in{\mathcal U}_n\}$. The family
$\{\{p_{D_n,U}\} : U\in{\mathcal U}_n\}$ is closure preserving. So,
since $X$ is T$_1$, the set $F_n = \{p_{D_n,U} : U\in{\mathcal
U}_n\}$ is discrete.
\end{proof}

\vspace{3mm} In particular, every metrizable space (or, more
generally, every T$_1$ M1-space (= a space with a $\sigma$-closure
preserving base) is DH$^+$-separable. But below we will see more: every M3 (= stratifiable) space is DH$^+$-separable.

Shapirovskii showed
\cite{Shap} that every space with a $\sigma$-point finite base is
d-separable.

\begin{question}
{\rm Is every space with a $\sigma$-point finite base D-separable? If yes, how about other properties in
Definitions~\ref{DefinitionOfDsep} and \ref{DefinitionOfDPlusSep}? }
\end{question}

\begin{proposition}
Let $X$ be a collectionwise Hausdorff discretely generated space with a $\sigma$-closed
discrete dense set. Then $X$ is DH$^+$-separable.
\end{proposition}

\begin{proof}
Let $H=\bigcup_{n\in\omega} H_n$ be dense in $X$ (where each $H_n$ is closed and discrete).
Without loss of generality we assume that $H_n \subset H_m$ whenever $n\leq m$. For every $n$,
fix a pairwise disjoint open expansion $\{U_{n,x} : x\in H_n\}$ of $H_n$.

At the $n$th inning ONE picks a dense $D_n\subset X$. Then Two, for every $x\in H_n$, picks a discrete
$F_{n,x} \subset D_n \cap U_{n,x}$ such that $x\in\overline{F_{n,x}}$ and sets $F_n = \bigcup_{x\in H_n} F_{n,x}$.
\end{proof}

\begin{corollary} \label{MNsigma}
Every monotonically normal $\sigma$-space is DH$^+$-separable.
\end{corollary}

\begin{proof}
Every $\sigma$-space has a $\sigma$-closed discrete dense set,
and every monotonically normal space is both collectionwise Hausdorff
(see \cite{GaryHandbook}) and discretely generated \cite{DTTW}.
\end{proof}
\begin{corollary}
Every stratifiable space is DH$^+$-separable.
\end{corollary}

\begin{proof}
Because a stratifiable space is both monotonically normal and a $\sigma$-space.
\end{proof}



Boaz Tsaban asked us in private communication whether a separable D-separable space has to be selectively separable. This can be disproved by taking the space $Seq(\mathcal{F})$ where $\mathcal{F}$ is any ultrafilter on $\omega$. Indeed, this space is countable, and hence it is trivially a $\sigma$-space. Moreover, it is monotonically normal by Theorem 3.2 of \cite{JSSresolvability}. If $\mathcal{F}$ is a Ramsey ultrafilter (which exists, for example, if one assumes CH), then $Seq(\mathcal{F})$ is even a topological group (see \cite{VaughanSeqP}). So we arrive to the following theorem:


\begin{theorem} \label{answerboaz}
There is a countable DH$^+$-separable space $X$ which is not selectively separable. Under CH the space $X$ can even be taken to be a topological group.
\end{theorem}

Yet the following is still unknown.

\begin{question}
{\rm Is there a compact separable $D$-separable non-selectively
separable space?}
\end{question}

\vspace{3mm} Monotone normality alone does not imply
D-separability. Indeed, it suffices to
consider a Suslin Line $\mathbb{L}$. $\mathbb{L}$ is monotonically normal
because it is linearly ordered, and it cannot even have a
$\sigma$-discrete dense set because every discrete set in $\mathbb{L}$ is
countable, but $\mathbb{L}$ is not separable. Moreover, it is easy to see that
for linearly ordered spaces, the three properties: d-separability,
D-separability and having a $\sigma$-discrete $\pi$-base, are
equivalent. This motivates the following question:

\begin{question}
{\rm Is it true that a monotonically normal space is $D$-separable
if and only if it is $d$-separable?}
\end{question}

We conclude the section with a partial positive result. The
principal tool in the proof of it is a theorem by Gartside stating that $\pi w(X) =
d(X) = hd(X)$ for $X$ having a monotonically normal compactification
\cite{Gartside}


\begin{theorem}
Suppose a space $X$ has a  monotonically normal compactification.
Then the following conditions are equivalent:
\begin{enumerate}
\item $X$ is d-separable;
\item $X$ is D-separable;
\item $X$ has a $\sigma$-disjoint $\pi$-base.
\end{enumerate}
\end{theorem}

\begin{proof}
Of course it is enough to prove (1) $\Rightarrow$ (3).
Call a non-empty open set {\it poor} if for every non-empty open
$V\subset U$, $\pi w(V) = \pi w(U)$. It is clear that every
non-empty open set contains a poor set and thus in every topological
space $X$ one can find a pairwise disjoint family of poor sets
$\mathcal U$ such that $X=\overline{\bigcup{\mathcal U}}$. Note that
$X$ is d-separable, or is D-separable, or is a space with a $\sigma$-disjoint $\pi$-base iff so is every element of $\mathcal U$. Moreover, if X has a monotonically normal compactification, then also every element of $\mathcal{U}$ does. So, without loss of generality
we can assume that $X$ itself is poor.

So let $X$ be a poor space with a monotonically normal
compactification and monotone normality operator $H$, and let $D=\bigcup_{n\in\omega} D_n$ be dense in $X$
where each $D_n$ is discrete. Without loss of generality, we can assume
that the sets $D_n$ are pairwise disjoint. Let $\mathcal P$ be a
$\pi$-base of $X$ of cardinality $\pi w(X)=\kappa$. By Gartside's
theorem (applied to the subspace $U$) and poverty of $X$, $|U\cap
D|=\kappa$ for every $U\in \mathcal P$. Then, enumerating $\mathcal
P$ and $D$, one easily defines an injection $f:{\mathcal P} \to D$
such that $f(U)\in U$ for every $U\in \mathcal P$. For $n\in\omega$,
put ${\mathcal P}_n = \{U\in{\mathcal P} : f(U)\in D_n\}$ and let $\{V(f(U)): U \in \mathcal{P}_n\}$ be a family of open sets such that $V(f(U)) \cap D_n=\{f(U)\}$ and $V(f(U)) \subset U$ for every $U \in \mathcal{P}_n$. For $U\in {\mathcal P}_n$, put $U^\prime = H(f(U),V(f(U))$. By the properties of the monotone normality operator $H$, the family $\mathcal{P}'_n=\{U': U \in \mathcal{P}_n \}$ is pairwise disjoint and hence $\bigcup_{n<\omega} \mathcal{P}'_n$ is a $\sigma$-disjoint $\pi$-base for $X$.
\end{proof}

\subsection{Subspaces, unions}

The following is straightforward:

\begin{proposition}\label{TrivialPropertiesOfDSpaces}
(1) Every open subspace, as well as every dense subspace, of a space
with one of the properties from Definitions~\ref{DefinitionOfDsep}
and \ref{DefinitionOfDPlusSep} has the same property.

(2) If $X$ has an open dense subspace with one of the properties
from Definitions~\ref{DefinitionOfDsep} and
\ref{DefinitionOfDPlusSep} then $X$ has the same property.

(3) A discrete sum of spaces with one of the properties from
Definitions~\ref{DefinitionOfDsep} and \ref{DefinitionOfDPlusSep}
 has the same property.
\end{proposition}

It was shown in \cite{GruenhageSakai} that selective separability,
R-separability and GN-separability are preserved by finite unions
(to see that this is not immediate, it might be enough to mention
that the question about H-separability remains open, and that SS$^+$
is not finitely additive \cite{BermanDowSecond}).

\begin{proposition}\label{PropLocallyFiniteSum}
A locally finite union of D-separable spaces is D-separable.
\end{proposition}

First with prove that this is the case for finite unions. The proof is a modification of the proof that selective separability is preserved by finite unions (see \cite{GruenhageSakai}).

\begin{lemma}\label{TheUnionOfTwo}
The union of two D-separable  spaces is D-separable.
\end{lemma}

\begin{proof}
Let $X=A \cup B$ where $A$ and $B$ are $D$-separable
and let $\{D_n: n \in \omega \} \subset X$ be a sequence of dense
sets. Let $U_n=X \setminus (\overline{(\bigcup_{i \geq n} D_i) \cap
A} \cup \overline{ \bigcup_{j < n} U_j})$ and $U=\bigcup_{n \in
\omega} U_n$. Then the sets $U_n$ are open in $X$ and pairwise
disjoint.

\medskip

\noindent {\bf Claim 1.}
For each $i \geq n$ the set $D_i \cap B \cap U_n$ is dense in $U_n$.

\begin{proof}[Proof of Claim 1] Because $U_n$ is
open, $D_i$ is dense in $X$, $X=A\cup B$, and for $i\geq n$, $D_i
\cap U_n \cap A = \varnothing$.
\renewcommand{\qedsymbol}{$\triangle$}
\end{proof}

\noindent {\bf Claim 2.}
There are discrete $G_n \subset D_n$ such that $\bigcup_{n \in
\omega} G_n$ is dense in $U$.

\begin{proof}[Proof of Claim 2] By
Proposition~\ref{TrivialPropertiesOfDSpaces}, part 1, the subspace
$B \cap U_n$ is D-separable for every $n \in \omega$. Therefore
there are discrete $G^n_i \subset D_i \cap B \cap U_n$ for every $i
\geq n$ such that $\bigcup_{i \geq n} G^n_i$ is dense in $B \cap
U_n$ and hence in $U_n$ (because $B\cap U_n$ is dense in $U_n$). Let
now $G_i=\bigcup_{n \leq i} G^n_i$. Since $G^n_i \subset U_n$ for
every $n \leq i$ and $\{U_n: n \leq i \}$ is a family of pairwise
disjoint open sets we have that $G_i$ is a discrete subset of $D_i$.
Moreover $\bigcup_{i \in \omega} G_i$ is dense in $U_n$ for every $n
\in \omega$, and hence in $U$.
\renewcommand{\qedsymbol}{$\triangle$}
\end{proof}

Let now $V=X \setminus \overline{U}$. We claim that $(\bigcup_{i
\geq n} D_i) \cap A$ is dense in $V$, and hence also in $A \cap V$.
Indeed, if $x \in V$, then $x \notin \overline{U}$, so $x \notin
U_n$ and $x \notin \overline{\bigcup_{j < n} U_j}$ for every $n \in
\omega$, which together imply that $x \in \overline{(\bigcup_{i \geq
n} D_i) \cap A}$ for every $n \in \omega$.

Now $A \cap V$ is D-separable, so there are discrete $H_n \subset
(\bigcup_{i \geq n} D_i) \cap A$ so that $\bigcup_{n \in \omega}
H_n$ is dense in $V \cap A$ and hence in $V$.

For each $x \in H_n$ let $i_n(x) \in \omega \setminus n$ be such
that $x \in D_{i_n(x)}$. Let $K_i=\{x: \exists n \in \omega (x \in
H_n$ and $i_n(x)=i) \}$. Then $K_i$ is a discrete subset of $D_i$
and $\bigcup_{i \in \omega} K_i=\bigcup_{n \in \omega} H_n$, hence
it is dense in $V$. Thus, if $G_n$ is as in Claim 2, then
$\bigcup_{n \in \omega} (G_n \cup K_n)$ is dense in $X$ and each
$G_n \cup K_n$ is discrete since $G_n \subset U$, $K_n \subset X
\setminus \overline{U}$ and $G_n$ and $K_n$ are both discrete. So
$X$ is D-separable.
\end{proof}

\begin{proof}[Proof of Proposition~\ref{PropLocallyFiniteSum}:] First of
all, by induction, Lemma~\ref{TheUnionOfTwo} can be extended to any
finite union.

Now, let $X=\bigcup{\mathcal Y}$ be a locally finite union, and let
each $Y\in \mathcal Y$ be D-separable. For $n\geq 1$ put $X_n =
\{x\in X :$ there is a neighborhood $U$ of $x$ such that
$|\{Y\in{\mathcal Y} : Y\cap U \neq \varnothing\}|\leq n\}$. Then
the sets $X_n$ are open in $X$, and $X=\bigcup_{n\in\omega} X_n$.
Put $Z_1=X_1$. For $n>1$, put $ Z_n = X_n
\setminus\overline{X_{n-1}}$. Then the sets $Z_n$ are open in $X$,
pairwise disjoint, and $X=\overline{\bigcup_{n\geq 1} Z_n}$.
Further, each $Z_n$ is a discrete union $Z=\bigsqcup\{Z_{n,A} :
A\subset {\mathcal Y}, |A|=n\}$ where $Z_{n,A} = \{x\in Z_n : \mbox{
there is a neighborhood } U \mbox{ of } x \mbox{ such that } |\{Y\in
{\mathcal Y} : Y\cap U \neq \varnothing\} = A\}$. Finally,
$X=\overline{\bigsqcup\{Z_{n,A} : n\geq 1,\hspace{1mm} A\subset
{\mathcal Y},\hspace{1mm} |A|=n\}}$ where the sets $Z_{n,A}$ are
open in $X$ and pairwise disjoint. By Lemma~\ref{TheUnionOfTwo}
(extended to arbitrary finite unions) each $Z_{n,A}$ is D-separable;
hence by Proposition~\ref{TrivialPropertiesOfDSpaces}, part 3, so is
$\bigsqcup\{Z_{n,A} : n\geq 1, A\subset {\mathcal Y}, |A|=n\}$,
hence by part 2 of the same proposition so is $X$.
\end{proof}

\begin{question}
{\rm Is every (locally) finite union of $DH$-separable spaces again
$DH$-separable?}
\end{question}

\subsection{Products}

The following two beautiful results witness how well-behaved d-separability is with respect to products.
\begin{theorem}\label{AVProducts}
{\rm (Arhangelskii, \cite{AVdSeparability})} Any product of
d-separable spaces is d-separable.
\end{theorem}

\begin{theorem}\label{JSzPower}
{\rm (Juh\'asz and Szentmikl\'ossy \cite{JSz})} For every T$_1$-space
$X$, $X^{d(X)}$ is d-separable.
\end{theorem}

Theorems \ref{AVProducts} and \ref{JSzPower} imply two obvious
corollaries:

\begin{corollary}\label{TwoCorollariesOnProducts}
(1) For every T$_1$-space $X$ there is $\kappa(X)$ such that for
every $\kappa\geq \kappa(X)$, $X^\kappa$ is d-separable.

(2) For every T$_1$-space $X$ there is a T$_1$-space $Z$ such that
the product $X\times Z$ is d-separable. If $X$ is Tychonoff, then so
is $Z$.
\end{corollary}

It is natural to ask if D-separability is (finitely or infinitely)
productive and if the analogue of Corollary~\ref{TwoCorollariesOnProducts}
is true for D-separability. It turns out that the product of two
D-separable spaces does not have to be D-separable, and that for
part 2 of Corollary~\ref{TwoCorollariesOnProducts} the answer is
affirmative while for part 1 the situation is almost the opposite.

\begin{theorem}\label{HighPowersAreNotDSeparable}
For every Tychonoff space $X$ with $|X|>1$ and every $\kappa$ there
is $\kappa^\prime \geq \kappa$ such that $X^{\kappa^\prime}$ is not
D-separable.
\end{theorem}

\begin{theorem}\label{SantinoOnProducts}
For every space $X$ there is a Tychonoff space $Z$ such that
$X\times Z$ is DH$^+$-separable.
\end{theorem}

But before proving the above theorems let's examine the case of
finite products.

\begin{theorem} \label{ProductSigmaDisjoint}
Let $X$ be a $D$-separable space (or has another property from
Definitions~\ref{DefinitionOfDsep} and \ref{DefinitionOfDPlusSep})
and $Y$ be a space having a $\sigma$-disjoint $\pi$-base. Then $X
\times Y$ is $D$-separable (or has the corresponding property).
\end{theorem}

\begin{proof} (For D-separability) Let
$\mathcal{B}=\bigcup_{n<\omega} \mathcal{B}_n$ be a
$\sigma$-disjoint $\pi$-base for $Y$ and let $\{D_k: k<\omega\}$ be
a countable sequence of dense subsets of $X\times Y$. Let
$\{B^n_\alpha: \alpha<\tau_n\}$ enumerate $\mathcal{B}_n$ and
$\{A_n: n < \omega \}$ be a partition of $\omega$. Observe that the
set $\pi_X(D_k \cap \pi^{-1}_Y(B^n_\alpha))$ is dense in $X$ for
every $k \in A_n$ and for every $\alpha \in \tau_n$. Fix $\alpha <
\tau_n$. Then for every $k \in A_n$ we can find a discrete set
$E^\alpha_k \subset \pi_X(D_k \cap \pi^{-1}_Y(B^n_\alpha))$ such
that $\bigcup_{k \in A_n} E^\alpha_k$ is dense in $X$. For every $x
\in E^\alpha_k$ pick a point $f(x) \in \pi_X^{-1}(x)$ such that
$\pi_Y(f(x)) \in B^n_\alpha$ and set $F^\alpha_k=\{f(x): x \in
E^\alpha_k \}$ and $F_k=\bigcup_{\alpha<\tau_k} F^\alpha_k$. Then
$F_k$ is discrete. In fact, let $(x,y) \in F_k$. Then $(x,y) \in
F^\alpha_k$ for some $\alpha < \tau_k$ so $(x,y)=f(x)$ for some $y
\in B^n_\alpha$. Moreover $x \in E^\alpha_k$. Now, since
$E^\alpha_k$ is discrete in $X$ there is an open $V \subset X$ such
that $V \cap E^\alpha_k=\{x\}$. Finally, observe that, since
$\mathcal{B}^n_\alpha$ is a disjoint family $(V \times B^n_\alpha)
\cap F_k=\{(x,y) \}$, which shows that $F_k$ is discrete. Moreover
$F_k \subset D_k$ and $\bigcup_{k<\omega} F_k$ is dense in $X \times
Y$, which proves that $X \times Y$ is $D$-separable.

The proofs for the other properties differ only by minor changes.
\end{proof}

\begin{example}
{\rm [CH]} The product of two countable selectively separable spaces
need not be D-separable.
\end{example}

{\it Proof:} Let $X$ be a selectively separable countable maximal
regular crowded space such that $X^2$ has no dense slim set, see
\cite{GruenhageSakai}. Let us check that the proof from
\cite{GruenhageSakai} that $X^2$ is not selectively separable
provides more: that $X^2$ is not D-separable. Enumerate $X=\{x_i :
i\in\omega\}$. For every $n \in \omega$ let $D_n=\{(x,y): x,y \notin
\{x_i: i \leq n \}\}$. Then $\{D_i: i \in \omega \}$ is a sequence
of dense sets in $X$. Let $E_n \subset D_n$ be a discrete set. Then
$\bigcup_{n \in \omega} E_n$ meets every cross-section in a finite
union of discrete sets. Now, in a crowded space, every discrete set
is nowhere dense and finite unions of nowhere dense sets are nowhere
dense. Therefore $\bigcup_{n \in \omega} E_n$ cannot be dense in
$X^2$, which proves that $X^2$ is not D-separable. $\Box$

\vspace{3mm} As a byproduct we get that under CH there exists a
countable non-D-separable space. However, one can construct such an
example even in ZFC.

\begin{example}\label{CountableNonDSeparable}
There is a dense countable subset $X\subset 2^{\mathfrak c}$ such that
$X$ is not D-separable.
\end{example}

\begin{lemma}\label{LemmaAboutInitialFace}
(1) For every countable subset $S\subset 2^{\mathfrak c}$, there is
$\alpha<\mathfrak c$ such that $\pi_{[0,\alpha)}|_S$ is a bijection.

(2) If a countable subset $S\subset 2^{\mathfrak c}$ is
$\sigma$-discrete, then this can be witnessed by a projection to
some initial face in $2^{\mathfrak c}$. That is, if
$S=\bigcup_{n\in\omega}S_n$ where each $S_n$ is countable and
discrete, then there is $\alpha<\mathfrak c$ such that
$\pi_{[0,\alpha)}(S_n)$ is discrete for each $n$, and
$\pi_{[0,\alpha)}|_S$ is injective.
\end{lemma}

\begin{proof}[Proof of lemma] (1) Pick countably many standard neighborhoods
of points of $S$ separating points of $S$ and use the fact that
$cf({\mathfrak c})>\omega$.

(2) Pick standard neighborhoods of points of $S$ witnessing
$\sigma$-discreteness.
\end{proof}

\vspace{3mm} {\it Construction of
Example~\ref{CountableNonDSeparable}:} First, it is easy to
construct pairwise disjoint dense countable subspaces $Y_n$,
$n\in\omega$  in $2^{\mathfrak c}$ such that for every two distinct
$y_1, y_2 \in Y=\bigcup_{n\in\omega}Y_n$, the set $I_{y_1, y_2} =
\{\alpha < {\mathfrak c} : y_1(\alpha) \neq y_2(\alpha)\}$ has
cardinality $\mathfrak c$. Using this, one can partition $\mathfrak c$ as
${\mathfrak c} = \cup\{C_A : A\subset Y\}$ so that each $C_A$ has
cardinality $\mathfrak c$, and for every two distinct $y_1, y_2 \in Y$,
$C_{\{y_1, y_2\}} \subset I_{y_1, y_2}$ (*).

Next, by induction on $0\leq \alpha<\mathfrak c$, we will construct
countable subspaces $Z_\alpha\subset 2^{\mathfrak c}$ that will take the
form $Z_\alpha = \{y_\alpha : y\in Y\}$. We also denote $Z_{\alpha,
n} = \{y_\alpha : y\in Y_n\}$, so $Z_\alpha = \bigcup_{n\in\omega}
Z_{\alpha,n}$. The points of $Z_\alpha$s are going to have the
following property: if $0\leq \gamma\leq\alpha \leq \beta < \mathfrak c$
 then for all $y\in Y$, $y_\beta(\gamma) =
y_\alpha(\gamma)$ (**).

To start the induction, we set $Z_0 = Y$, that is $y_0=y$ for all
$y\in Y$.

Now let $0<\alpha < \mathfrak c$, and suppose $X_\gamma$s have been
defined for all $\gamma<\alpha$. Let $y\in Y$. To define the
corresponding point $y_\alpha \in Z_\alpha$, we have to define
$y_\alpha(\gamma)$ for all $\gamma$, $0\leq \gamma < \mathfrak c$. If
$0\leq \gamma < \alpha$, then we set $y_\alpha(\gamma) =
y_\gamma(\gamma)$ (and thus condition (**) continues to hold.)

To define $y_\alpha(\alpha)$ we need some auxiliary notation. By the
previous, we have in fact defined $\pi_{[0,\alpha)}(y_\alpha)$ for
all $y\in Y$. For a subset $B\subset Y$, set $B_{<\alpha} =
\{\pi_{[0,\alpha)}(y_\alpha) : y\in B\} \subset 2^{[0,\alpha)}$. We
have $\alpha\in C_A$ for some $A\subset Y$. If all the following
conditions hold:
\begin{itemize}
\item (1) $A$ is infinite,
\item (2) the mapping $A\to A_{<\alpha}$ given by $y\mapsto \pi_{[0,\alpha)}(y_\alpha)$ is a bijection,
\item (3) for every $n\in\omega$, $(A\cap Y_n)_{<\alpha}$ is discrete,
\end{itemize}
then we set $y_\alpha(\alpha)=0$ for all $y\in A$. Otherwise we set
$y_\alpha(\alpha)=y(\alpha)$.

Finally, for all $\gamma$ with $\alpha < \gamma < \mathfrak c$, we set
$y_\alpha(\gamma) = y(\gamma)$. This concludes the construction of
$Z_\alpha$.

Now we define the countable subspace $X\subset 2^{\mathfrak c}$,
$X=\{\tilde{y} : y\in Y\}$ by setting $\tilde{y}(\alpha) =
y_\alpha(\alpha)$ for all $y\in Y$. It follows from (**) that
$\tilde{y}(\gamma) = y_\alpha(\gamma)$ whenever $0\leq \gamma \leq
\alpha < \mathfrak c$. For $n\in \omega$, we set $X_n =\{\tilde{y} :
y\in Y_n\}$, thus we have $X = \cup_{n\in\omega} X_n$.

\medskip

\noindent {\bf Claim 1.} The mapping $y \mapsto \tilde{y}$ from $Y$ onto $X$
is a bijection.

\begin{proof}[Proof of Claim 1] Indeed, if $y_1, y_2$ be distinct elements of $Y$,
then by our construction, since by (*) $C_{\{y_1, y_2\}} \subset
I_{y_1, y_2}$, we have $\tilde{y_1}(\alpha) = y_1(\alpha) \neq
y_2(\alpha) = \tilde{y_2}(\alpha)$ for every $\alpha\in C_{\{y_1,
y_2\}}$.
\renewcommand{\qedsymbol}{$\triangle$}
\end{proof}

\noindent {\bf Claim 2.} Each $X_n$ is dense in $2^{\mathfrak c}$ (and thus in
$X$).
\begin{proof}[Proof of Claim 2]
Indeed, let $F\subset \mathfrak c$ be finite, and let $\varphi \in
2^F$. We have to find $\tilde{y} \in X_n$ such that $\tilde{y}|_F =
\varphi$. For each $i\in F$, there is $A_i \subset Y$ such that
$i\in C_{A_i}$. Put ${\mathcal A} = \{A_i : i\in F$ and conditions
(1), (2), (3) were satisfied when the $i$th coordinates of the
points of $X_i$ were defined$\}$. Pick $\alpha^*$ with $\max(F) <
\alpha^* < \mathfrak c$. Using Lemma~\ref{LemmaAboutInitialFace},~(1),
we can assume that $\pi_{[0,\alpha^*)}|_X$ is a bijection. Then
$T=\pi_{[0,\alpha^*)}(\cup\{A_i\cap Y_n : A_i\in{\mathcal A}\})$ is
a finite union of discrete subspaces of $2^{[0,\alpha^*)}$, and thus
$T$ is nowhere dense in $2^{[0,\alpha^*)}$. So $T^\prime =
\pi_{[0,\alpha^*)}(Y_n) \setminus T$ is dense in $2^{[0,\alpha^*)}$.
Pick $t\in T^\prime$ with $t|_F = \varphi$ and $y\in Y_n$ with
$\pi_{[0,\alpha^*)}(y)=t$. Then $\tilde{y}\in X_n$, and
$\tilde{y}|_F = y|_F = t|_F = \varphi$.
\renewcommand{\qedsymbol}{$\triangle$}
\end{proof}

\noindent {\bf Claim 3.} For any choice of discrete $S_n \subset X_n$,
$n\in\omega$, the set $S=\cup_{n\in\omega} S_n$ is not dense in
$2^{\mathfrak c}$ (and thus not dense in $X$.)

\begin{proof}[Proof of Claim 3]
By Lemma~\ref{LemmaAboutInitialFace},~(2), there is $\alpha^*<\mathfrak c$
such that $\pi_{[0,\alpha^*)}(S_n)$ is discrete for each $n$, and
$\pi_{[0,\alpha^*)}|_S$ is injective. Put $A=\{y\in Y : \tilde{y}\in
S\}$. Pick $\alpha^{**} \in C_A$ so that  $\alpha^{**} \geq
\alpha^*$. Then $\tilde{y}(\alpha^{**}) = 0$ for every $\tilde{y}\in
S$, and thus $S$ is not dense in $2^{\mathfrak c}$.
\renewcommand{\qedsymbol}{$\triangle$}
\end{proof}

Claims 2 and 3 show that $X$ is as was desired. $\Box$

\begin{corollary}\label{2totheCIsNotDSeparable}
$2^{\mathfrak c}$ is not D-separable.
\end{corollary}

\begin{question}
{\rm What is $\mathfrak{cds}=\min\{\tau : 2^\tau$ contains a dense countable
subspace which is not D-separable$\}$?}
\end{question}

\begin{question}
{\rm What is $\mathfrak{ds}=\min\{\tau : 2^\tau$ is not D-separable$\}$?}
\end{question}

\begin{question}
{\rm Is $\mathfrak{cds}=\mathfrak{ds}$?}
\end{question}

\begin{question}
{\rm Is it true that for every separable Tychonoff space $X$ there
is a separable Tychonoff space $Y$ such that $X\times Y$ is
D-separable?}
\end{question}

\subsubsection{Proof of Theorem \ref{HighPowersAreNotDSeparable}  }

It suffices to show that (\dag) for every Tychonoff $X$ there is
$\tau$ such that $X^\tau$ is not D-separable. Indeed, applying
(\dag) to $X^\prime=X^\kappa$ we get the original statement of the
theorem. In the case of a finite $X$, $\tau=\mathfrak c$ works by an
easy modification of the argument from
Example~\ref{CountableNonDSeparable}, so we assume $\lambda=|X|$ is
infinite. Next, if $D$ is dense in $X$ and $D^\tau$ is not
D-separable, then neither is $X^\tau$. So we can pass from $X$ to a
dense subspace of minimal cardinality and thus assume $|X|=d(X)$
when proving the following

\begin{theorem} \label{ahighpowerisnotDseparable}
For every Tychonoff $X$, $X^{2^{d(X)}}$ is not D-separable.
\end{theorem}

{\it Proof.} The argument is parallel to one from
Example~\ref{CountableNonDSeparable}, so we will omit some details.
Fix a point $x_0\in X$.

Let $\tau = 2^\lambda$ (where $\lambda = |X| = d(X)$). Since ${\rm
cf}(\tau) > \lambda$ we get the following:

\begin{lemma}\label{LemmaAboutInitialFaceGen}
(1) For every subset $S\subset X^\tau$, such that $|S| \leq \lambda$ there is $\alpha<\tau$ such
that $\pi_{[0,\alpha)}|_S$ is a bijection.

(2) If a subset $S\subset 2^{\tau}$ such that $|S|\leq\lambda$ is
$\sigma$-discrete, then this can be witnessed by a projection to
some initial face in $X^\tau$. That is, if
$S=\bigcup_{n\in\omega}S_n$ where each $S_n$ is discrete, then there
is $\alpha<\tau$ such that $\pi_{[0,\alpha)}(S_n)$ is discrete for
each $n$, and $\pi_{[0,\alpha)}|_S$ is injective.
\end{lemma}

The routine proof of the next lemma is omitted.

\begin{lemma}
There exist pairwise disjoint dense subspaces $Y_n$, $n\in\omega$ in
$X^\tau$ such that $|Y_n| \leq \lambda$ and for every two distinct $y_1, y_2
\in Y=\bigcup_{n\in\omega}Y_n$, the set $I_{y_1, y_2} = \{\alpha <
\tau : y_1(\alpha) \neq y_2(\alpha)\}$ has cardinality $\tau$.
\end{lemma}

Using this, one can partition $\tau$ as ${\tau} = \bigcup\{C_A :
A\subset Y\}$ so that each $C_A$ has cardinality $\tau$, and for
every two distinct $y_1, y_2 \in Y$, $C_{\{y_1, y_2\}} \subset
I_{y_1, y_2}$ (*).

Next, by induction on $0\leq \alpha<\tau$, we will construct
$\lambda$-sized subspaces $Z_\alpha\subset X^\tau$ that will take the form
$Z_\alpha = \{y_\alpha : y\in Y\}$. We also denote $Z_{\alpha, n} =
\{y_\alpha : y\in Y_n\}$, so $Z_\alpha = \bigcup_{n\in\omega}
Z_{\alpha,n}$. The points of $Z_\alpha$s are going to have the
following property: if $0\leq \gamma\leq\alpha \leq \beta < \tau$
 then for all $y\in Y$, $y_\beta(\gamma) =
y_\alpha(\gamma)$ (**).

To start the induction, we set $Z_0 = Y$, that is $y_0=y$ for all
$y\in Y$.

Now let $0<\alpha < \tau$, and suppose $X_\gamma$s have been defined
for all $\gamma<\alpha$. Let $y\in Y$. To define the corresponding
point $y_\alpha \in Z_\alpha$, we have to define $y_\alpha(\gamma)$
for all $\gamma$, $0\leq \gamma < \tau$. If $0\leq \gamma < \alpha$,
then we set $y_\alpha(\gamma) = y_\gamma(\gamma)$ (and thus
condition (**) continues to hold.)

To define $y_\alpha(\alpha)$ we need some auxiliary notation. By the
previous, we have in fact defined $\pi_{[0,\alpha)}(y_\alpha)$ for
all $y\in Y$. For a subset $B\subset Y$, set $B_{<\alpha} =
\{\pi_{[0,\alpha)}(y_\alpha) : y\in B\} \subset X^{[0,\alpha)}$. We
have $\alpha\in C_A$ for some $A\subset Y$. If all the following
conditions hold:
\begin{itemize}
\item (1) $A$ is infinite,
\item (2) the mapping $A\to A_{<\alpha}$ given by $y\mapsto \pi_{[0,\alpha)}(y_\alpha)$ is a bijection,
\item (3) for every $n\in\omega$, $(A\cap Y_n)_{<\alpha}$ is discrete,
\end{itemize}
then we set $y_\alpha(\alpha)=x_0$ for all $y\in A$. Otherwise we
set $y_\alpha(\alpha)=y(\alpha)$.

Finally, for all $\gamma$ with $\alpha < \gamma < \tau$, we set
$y_\alpha(\gamma) = y(\gamma)$. This concludes the construction of
$Z_\alpha$.

So we have $Z_\alpha$ satisfying (**) for all $\alpha<\tau$. Now we
define the subspace $\tilde{Y}\subset X^\tau$ by
$\tilde{Y}=\{\tilde{y} : y\in Y\}$ where $\tilde{y}(\alpha) =
y_\alpha(\alpha)$ for all $y\in Y$. It follows from (**) that
$\tilde{y}(\gamma) = y_\alpha(\gamma)$ whenever $0\leq \gamma \leq
\alpha < \tau$. For $n\in \omega$, we set $\tilde{Y}_n =\{\tilde{y}
: y\in Y_n\}$, thus we have $\tilde{Y} = \cup_{n\in\omega}
\tilde{Y}_n$.
\medskip

\noindent {\bf Claim 1.} The mapping $y \mapsto \tilde{y}$ from $Y$ onto
$\tilde{Y}$ is a bijection.

 \begin{proof}[Proof of Claim 1]
Indeed, if $y_1, y_2$ be distinct
elements of $Y$, then by our construction, since by (*) $C_{\{y_1,
y_2\}} \subset I_{y_1, y_2}$, we have $\tilde{y_1}(\alpha) =
y_1(\alpha) \neq y_2(\alpha) = \tilde{y_2}(\alpha)$ for every
$\alpha\in C_{\{y_1, y_2\}}$.
\renewcommand{\qedsymbol}{$\triangle$}
\end{proof}

\noindent {\bf Claim 2.} Each $\tilde{Y}_n$ is dense in $X^\tau$ (and thus in
$\tilde{Y}$).

\begin{proof}[Proof of Claim 2]
Indeed, let $F\subset \tau$ be finite, and let
$\varphi \in ({\mathcal T}\setminus \{\varnothing\})^F$ (where
$\mathcal T$ is the topology of $X$). We have to find $\tilde{y} \in
\tilde{Y}_n$ such that (+) $\tilde{y}(i) \in \varphi(i)$ for every
$i\in F$. For each $i\in F$, there is $A_i \subset Y$ such that
$i\in C_{A_i}$. Put ${\mathcal A} = \{A_i : i\in F $ and conditions
(1), (2), (3) were satisfied when the $i$th coordinates of the
points of $X_i$ were defined$\}$. Pick $\alpha^*$ with $\max(F) <
\alpha^* < \tau$. Using Lemma~\ref{LemmaAboutInitialFaceGen},~(1),
we can assume that $\pi_{[0,\alpha^*)}|_{\tilde{Y}}$ is a bijection.
Then $T=\pi_{[0,\alpha^*)}(\cup\{A_i\cap Y_n : A_i\in{\mathcal
A}\})$ is a finite union of discrete subspaces of
$X^{[0,\alpha^*)}$, and thus $T$ is nowhere dense in
$X^{[0,\alpha^*)}$. So $T^\prime = \pi_{[0,\alpha^*)}(Y_n) \setminus
T$ is dense in $X^{[0,\alpha^*)}$. Pick $t\in T^\prime$ with $t(i)
\in \varphi(i)$ for every $i\in F$. There is  $y\in Y_n$ with
$\pi_{[0,\alpha^*)}(y)=t$. Then $\tilde{y}\in \tilde{Y}_n$, and
$\tilde{y}|_F = y|_F$, so $\tilde{y}$ satisfies (+).
\renewcommand{\qedsymbol}{$\triangle$}
\end{proof}

{\bf Claim 3.} For any choice of discrete $S_n \subset \tilde{Y}_n$,
$n\in\omega$, the set $S=\bigcup_{n\in\omega} S_n$ is not dense in
$X^\tau$ (and thus not dense in $X$.)
\begin{proof}[Proof of Claim 3]
By Lemma~\ref{LemmaAboutInitialFaceGen},~(2), there is $\alpha^*<\tau$
such that $\pi_{[0,\alpha^*)}(S_n)$ is discrete for each $n$, and
$\pi_{[0,\alpha^*)}|_S$ is injective. Put $A=\{y\in Y : \tilde{y}\in
S\}$. Pick $\alpha^{**} \in C_A$ so that  $\alpha^{**} \geq
\alpha^*$. Then $\tilde{y}(\alpha^{**}) = x_0$ for every
$\tilde{y}\in S$, and thus $S$ is not dense in $X^\tau$.
\renewcommand{\qedsymbol}{$\triangle$}
\end{proof}

Claims 2 and 3 show that $\tilde{Y}$ is not D-separable. Since
$\tilde{Y}$ is dense in $X^\tau$ it follows that $X^\tau$ is not
D-separable. $\Box$

\begin{question}
{\rm Is it true that for every Tychonoff space $X$ there is $\kappa$
such that for all $\kappa^\prime \geq \kappa$, $X^{\kappa^\prime}$
is not D-separable?}
\end{question}

\subsubsection{Proof of Theorem \ref{SantinoOnProducts}  }

 More specifically, we will prove:

\begin{theorem} \label{everyspacehasaDseparableproduct}
Let $X$ be any space, and let $Y$ be any space such that $\pi
w(X)\leq \pi w(Y)=\kappa$ and $Y$ contains a cellular family of size
$\kappa$. Then $X\times Y^\omega$ is DH$^+$-separable.
\end{theorem}

(Then, for Theorem~\ref{SantinoOnProducts}, one can take
$Z=Y^\omega$. As $Y$, one can take the discrete space of size $\pi
w(X)$ or a one-point compactification of such a space, so $Z$ in
Theorem~\ref{SantinoOnProducts} can be in addition assumed compact.)

\begin{proof}
Let $\mathcal U$ and $\mathcal V$ be $\pi$-bases of $X$
and $Y$ having minimal size. Let $\{C_\alpha : \alpha<\kappa\}$ be a
cellular family in $Y$. For $m\in\omega$, let $\{e^m_\alpha :
\alpha<\kappa\}$ be an enumeration of ${\mathcal U}\times {\mathcal
V}^m$.

On the $m$th move, One chooses a dense subspace $S_m\subset X\times
Y^\omega$, and Two, for every $m\in\omega$ and $\alpha<\kappa$
selects $d_\alpha^m \in (\prod_{n\in\omega} W_n^{\alpha,m})\cap S_m$
where $W_n^{\alpha,m} = e^m_\alpha(n)$ for $n\leq m$,
$W_{m+1}^{\alpha,m} = C_\alpha$ and $W_n^{\alpha,m} = Y$ for
$n>m+1$. Let $D_m = \{d_\alpha^m : \alpha<\kappa\}$. Then $D_m
\subset S_m$, $D_m$ is discrete, and $D_m$ intersects every
non-empty open set in $X\times Y^\omega$ that depends only on the
first $m+1$ coordinates. Thus every non-empty open set in $X\times
Y^\omega$ intersects all but finitely many $D_m$s.
\end{proof}

A consequence of the above Theorem is that there is no single
cardinal $\kappa$ such that $X^\kappa$ is \underline{not}
$D$-separable for every space $X$.

\begin{corollary}
For every $\kappa$, and every $\lambda\leq \kappa$,
$(D(\kappa))^\lambda$ is D-separable (where $D(\kappa)$ is the
discrete space of cardinality $\kappa$).
\end{corollary}

\begin{proof}
This is trivial if $\kappa$ is finite. So assume $\kappa$ is infinite. Now set $X=Y=D(\kappa)^\omega$ in Theorem $\ref{everyspacehasaDseparableproduct}$ and observe that $(D(\kappa)^\kappa)^\omega$ and $D(\kappa)^\kappa$ are homeomorphic.
\end{proof}


Another notable consequence of Theorem
$\ref{everyspacehasaDseparableproduct}$ is the fact that the
$\omega$-power of any linearly ordered space is $D$-separable. This
follows from the following result of Petr Simon.

\begin{lemma}
\cite{Simon} Let $X$ be a linearly ordered topological space. Then
$X^2$ contains a cellular  family of size $d(X)$.
\end{lemma}

\begin{corollary}
Let $X$ be a linearly ordered topological space. Then $X^\omega$ is
$DH^+$-separable.
\end{corollary}

\begin{proof}
Since $\pi w(X)=d(X)$ in linearly ordered spaces $X^2$ contains a pairwise disjoint open family of size $\pi w(X)$. Now let $Y=X^2$ in Theorem $\ref{everyspacehasaDseparableproduct}$.
\end{proof}

So, although a Suslin Line is not even $d$-separable, its
$\omega$-power is $DH^+$-separable.

\subsection{Some more open problems}

Tkachuk presented a large collection of sufficient conditions and
necessary conditions of d-separability of $C_p(X)$ in
\cite{Tka2005}. Tkachuk gave a CH example of a compact space $X$
with a non-d-separable $C_p(X)$ and  asked for a ZFC example of a
Tychonoff space or even a compact space $X$ with non-d-separable
$C_p(X)$. A Tychonoff ZFC example was presented in \cite{JSz}.

\begin{problem}
{\rm (1) Characterize $X$ such that $C_p(X)$ is D-separable.

(2) More specifically, suppose $C_p(X)$ is d-separable. Under
what additional conditions on $X$ is $C_p(X)$ D-separable? }
\end{problem}

Recall that a compact space $X$ is selectively separable iff $X$ has
a countable $\pi$-base (\cite{BBMT}). This is a consequence of the fact that
a compact space $X$ has a countable $\pi$-base iff every dense
subspace of $X$ is separable \cite{JShelah}.

Let $dd(Y)$ be the least cardinal $\kappa$ such that $Y$ has a dense set which is the union of $\kappa$ many discrete sets. Let $d\delta(X)=\sup \{dd(D): D$ is dense in $X \}$. Let $d\pi(X)$ be the least cardinal $\kappa$ such that $X$ has a $\pi$-base which is the union of $\kappa$ many disjoint collections.

\begin{conjecture} \ \label{compactconjectures}
\begin{enumerate}
\item \label{weakconjecture} A compact space $X$ is D-separable iff $X$ has a $\sigma$-disjoint
$\pi$-base.
\item \label{strongconjecture} Let $X$ be a compact space. Then $d\delta(X)=d\pi(X)$.
\end{enumerate}
\end{conjecture}

By Theorem \ref{ProductSigmaDisjoint} if Conjecture $\ref{compactconjectures}$, $\ref{weakconjecture}$ is true then
the answer to the following question is positive.

\begin{question}
{\rm Is the product of two compact D-separable spaces still
D-separable?}
\end{question}

Recall that a space is called an $L$-space if it is hereditarily Lindel\"of but not separable. Tkachuk \cite{Tka2005} constructed under CH an L-space $X$ such that $X^2$ is $d$-separable. Later on, Moore \cite{Justin} showed that a slight modification of his ZFC example of an L-space provides a ZFC example of an $L$-space with a $d$-separable square.

\begin{question}
{\rm Is there a non-D-separable space $X$ such that $X^2$ is
D-separable? Is there even a non-d-separable space with this
property?}
\end{question}

\begin{question}
{\rm Is there (in any model of ZFC) an example of an L-space with a
D-separable square?}
\end{question}

Note that replacing D-separability with selective separability both questions have easily a negative answer.


Also, the influence of convergence properties on D-separability is not clear yet.


\begin{question}\ {\rm
\begin{enumerate}
\item Is every Fr\'echet $d$-separable space $D$-separable? What about $\Sigma(2^\mathfrak{\kappa})$?
\item Is there a sequential $d$-separable (separable, countable) non-D-separable space?
\item Is there a Whyburn $d$-separable (separable, countable) non-D-separable space?
\end{enumerate}}
\end{question}

\section{More on maximal (and submaximal) spaces}\label{MiscellaniaSEction}

We conclude with some remarks on the interesting case of maximal and submaximal spaces. In a submaximal space every dense set is open, so in some sense
dense sets are \lq\lq big". This  implies  \lq\lq a lot of freedom"
in choosing a finite set and this in turn could suggest that a
maximal space can easily be selectively separable, but we will see
below that often things go differently.

In \cite{BBMT} it was shown that assuming ${\mathfrak d}=\omega_1$ there
is a maximal regular space which is not selectively separable, and
it was asked (1) whether or not such an example is possible within
ZFC, and (2) is it true (at least consistently) that every countable
maximal regular space is not selectively separable? Here is the
progress obtained since then:

\begin{theorem} \ \label{ProgressOnMaximal}
\begin{enumerate}
\item {\rm (Barman and Dow, \cite{BermanDowTopProc})} There is
(within ZFC) a countable maximal regular space which is not
selectively separable.

\item {\rm (Barman and Dow, \cite{BermanDowTopProc}, Repov\v{s} and Zdomskyy,
\cite{RepovsZdomskyy})} Consistently, there is no submaximal  SS
space (specifically, the existence of such a space implies the
existence of a separable P-set in $\omega^*$ while the existence of
a ccc P-set in $\omega^*$ is known to be independent from ZFC from
\cite{FrankiewiczShelahZbierskiPsets}).\footnote{The result is
stated in \cite{BermanDowTopProc}, \cite{RepovsZdomskyy} only for
maximal spaces, but it is easy to notice that the argument uses only
submaximality.}

\item {\rm (Barman and Dow, \cite{BermanDowTopProc})} {\rm
[MA$_{ctble}$]} There exists a maximal regular countable selectively
separable space. (So the existence of a maximal regular selectively
separable space is independent of ZFC.)

\item {\rm (Gruenhage and Sakai, \cite{GruenhageSakai})} {\rm [CH]}
There is a maximal space $X$ such that $X$ is R-separable but $X^2$
is not selectively separable.\footnote{The first author and
Gruenhage obtained the same result under a weaker assumption
MA$_{ctble}$.}

\item {\rm (Barman and Dow, \cite{BermanDowTopProc})} Every crowded
SS$^+$ space is resolvable. Hence no maximal space is
SS$^+$.\footnote{And one can see from the argument in
\cite{BermanDowTopProc} that,  more generally, no crowded submaximal
space) is SS$^+$.}
\end{enumerate}
\end{theorem}

Many results on maximal regular spaces, in particular the
construction of Barman and Dow (Theorem~\ref{ProgressOnMaximal},
part~1 above) are based on the following tool found by van~Douwen:

\begin{theorem} {\rm \cite{vDMaximal}}
For any countable regular crowded space $(X,\tau)$ there is a
stronger regular topology $\sigma$ such that the space $(X,\sigma)$
has a dense subset which is  a maximal space.
\end{theorem}

Below we present an alternative proof of
Theorem~\ref{ProgressOnMaximal}, part~1 (based on the space
$Seq(\mathcal F)$ as the starting point), and construct a maximal
regular countable SS space using a weaker assumption than in
Theorem~\ref{ProgressOnMaximal}, part~3, namely ${\mathfrak d}={\mathfrak
c}$. Then we discuss maximal D-separable spaces.

\begin{theorem}
There exists a countable regular maximal space which
is not selectively separable.
\end{theorem}

\begin{proof}
Start by letting $(X,\tau)=Seq(\mathcal F)$ and   fix
the sequence of dense subsets $\{H_n:n<\omega\}$, where
$H_n=\bigcup\{{}^k\omega:n\le k <\omega\}$.

\noindent {\it Step 1}:  use van Douwen's theorem to find
$\sigma\supset \tau$ and a dense subset $Z$ of $(X,\sigma)$ which is
a regular maximal space.

\noindent  {\it Step 2}:  since  each $H_n$ has a closed scattered
complement in $(X,\tau)$, it follows that $H_n$ remains dense  and
open in $(X,\sigma)$ and so  the set $D_n=Z\cap H_n$ is dense in
$Z$.

\noindent {\it Step 3}:  the sequence $\{D_n:n<\omega\}$ cannot have
a \lq\lq good selection" because it would be also a \lq\lq good
selection" for the sequence $\{H_n:n<\omega\}$ in $(X,\tau)$.
\end{proof}

\vspace{3mm} Recall that first Gruenhage under [CH] (later included
in \cite{GruenhageSakai}) and then Barman and Dow under $MA_{ctble}$
have shown the existence of a countable regular maximal selectively
separable space. Gruenhage's construction gives a stronger result: a
maximal R-separable space  whose square is not selectively
separable. We are going to show that, with respect  to the weaker
task to have just a maximal selectively separable space, $\mathfrak
d=\mathfrak c$ suffices.  The construction we present below follows the
pattern of that of Gruenhage.

\begin{lemma}\label{AngelosLemma1AboutMaximal}
Let $X$ be a space and $x\in X$. If $t(x,X)=\omega$
and $\chi(x,X)<\mathfrak d$, then $X$ has countable fan tightness at
$x$.
\end{lemma}

\begin{proof}
Let $\{A_n:n<\omega\}$ be a sequence of sets such that
$x\in \overline {A_n}$ for each  $n$. Since $t(x,X)=\omega$, we may
assume  each $A_n$ countable and write $A_n=\{a_{n,k}:k<\omega\}$.
Let $\{U_\alpha :\alpha <\kappa \}$ be a local base at $x$ with
$\kappa <\mathfrak d$. For any $\alpha $ we may define a function
$f_\alpha \in{}^\omega\omega$ by letting $f_\alpha (n)=\min\{k:
a_{n,k}\in U_\alpha \cap A_n\}$.  Since $\kappa <\mathfrak d$, the
family $\{f_\alpha :\alpha <\kappa \}$ cannot be dominating and so
there exists  $g\in {}^\omega\omega$ such that  the set $\{n:
f_\alpha(n)\le g(n)\}$ is infinite for each $\alpha $. Now, by
letting  $F_n=\{a_{n,k}:k\le g(n)\}$,  we may easily check that
$x\in \overline {\bigcup\{F_n:n<\omega\}}$.
\end{proof}

\vspace{3mm} In \cite{BellaMalykhinTightnessAndResolvability} it is
shown that any crowded space of countable fan tightness is
$\omega$-resolvable. So we have:

\begin{corollary}\label{AngelosCorollary1AboutMaximal}
A countable crowded space of weight less than
$\mathfrak d$ is $\omega$-resolvable.
\end{corollary}

The above corollary is the main ingredient in the proof of the
following:

\begin{lemma}\label{AngelosLemma2AboutMaximal} {\rm [${\mathfrak d}={\mathfrak c}$]}
Let $(X,\tau)$ be a countable crowded regular space of weight
$\leq\kappa$ where $\kappa < \mathfrak c$. Then:

(1) If $A$ is a dense subset of $(X,\tau)$, then there is an
enlargement $\sigma_1$ of $\tau$ such that $(X,\sigma_1)$ is a
regular crowded space of weight $\leq\kappa $, $A\in \sigma_1$ and
each dense open set in $(X,\tau)$ remains dense in $(X,\sigma_1)$;

(2) If $A$ is a crowded subset in $(X,\tau)$, then there exists an
enlargement $\sigma_2$ of $\tau$ such that $(X,\sigma_2)$ is a
regular crowded space of weight $\leq\kappa $, $A$ is either open in
$(X,\sigma_2)$ or it has an isolated point in $(X,\sigma_2)$ and
each dense open set in $(X,\tau)$ remains dense in $(X,\sigma_2)$.
\end{lemma}

\begin{proof}
Part 1: by Corollary~\ref{AngelosCorollary1AboutMaximal} the subspace $A$ is
$\omega$-resolvable and  we may write $A=\bigcup\{A_n:n<\omega\}$,
where each $A_n$ is dense and $A_i\cap A_j=\varnothing$ whenever
$i\ne j$. $\sigma_1$ is the topology on $X$ generated by
$\tau\cup\{A_n:n<\omega\}\cup\{X\setminus A_n:n<\omega\}$.

Part 2: if $A$ is dense in $X$, then we may argue as in part 1. If
not, let $V={\rm Int}(X\setminus A)$ and  consider the topology
$\tau'$ generated by $\tau\cup\{\overline V\}$. $\tau'$ is regular
and any dense open set  in $\tau$ remains dense in $\tau'$. If
$\overline V\cap A=\varnothing$, then apply part 1 to the space
$(X,\tau')$ and the dense set  $\overline V\cup A$. In the resulting
topology $\sigma _2$ the set $A$ is open. If $\overline V\cap A\ne
\varnothing$, then pick a point $p\in \overline V\cap A$ and apply
part 1 to the space $(X,\tau')$ and the dense set $X\setminus
(\overline V\cap A\setminus \{p\})$. In the resulting topology
$\sigma_2$ the point $p$ is isolated in $A$.
\end{proof}

\begin{theorem}\label{AngelosTheorem4AboutMax} {\rm [${\mathfrak d}=\mathfrak c$]}
There exists a countable
regular maximal selectively separable space.
\end{theorem}

\begin{proof}
Let $\tau_0$ be a regular crowded second countable
topology on the set $\omega$. List all infinite subsets of $\omega$
as $\{A_\alpha :\alpha <\mathfrak c\}$ and all $\omega$-sequences of
subsets of $\omega$ as $\{\langle D^\alpha _n:n<\omega\rangle
:\alpha<\mathfrak c\}$ (in the latter each element is listed $\mathfrak
c$-many times). For any $\alpha <\mathfrak c$ we will construct a
crowded regular topology $\tau_\alpha $ on $\omega$ in such a way
that:
\begin{enumerate}
\item if $\beta<\alpha $ then $\tau_\beta\subseteq \tau_\alpha $
and any dense open set in $\tau_\beta$ remains dense in $\tau_\alpha
$;
\item  the weight of $\tau_\alpha $ is at most $|\alpha
|+\omega$;
\item  if $A _\alpha $ is dense in $\tau_\alpha$,  then
$A_\alpha $   is dense open in $\tau_{\alpha +1}$;
\item if $A_\alpha $ is not dense but crowded in
$\tau_\alpha $, then either $A_\alpha $ is open in $\tau_{\alpha
+1}$ or $A_\alpha $ has an isolated point if $\tau_{\alpha +1}$;

\item if $\langle D_n^\alpha :n<\omega\rangle$  is a sequence of
dense open sets in $\tau_\alpha $, then   there are finite sets
$F_n^\alpha \subseteq D_n^\alpha $ such that the set
$\bigcup\{F_n^\alpha :n<\omega\}$  is dense open in $\tau_{\alpha
+1}$.
\end{enumerate}
Suppose to have already defined  topologies $\tau_\beta$ and
sequences $\langle F_n^\beta:n<\omega\rangle$ for $\beta<\alpha $
satisfying the above conditions. If $\alpha $ is a limit ordinal,
then we take as $\tau_\alpha $ the topology generated by
$\bigcup\{\tau_\beta:\beta<\alpha \}$. In this case, only condition
2 needs to be checked. Now, assume $\alpha =\gamma+1$. If $A_\gamma$
is crowded, then apply Lemma~\ref{AngelosLemma2AboutMaximal}  to get
a topology $\tau'$ ($\tau'$ is either $\sigma_1$ or $\sigma_2$ from
Lemma~\ref{AngelosLemma2AboutMaximal} according to the fact that
$A_\gamma $ is or is not dense in $\tau_\gamma$). Next, if $\langle
D_n^\gamma:n<\omega\rangle$ is a sequence of dense open sets in
$\tau_\gamma$ (and so even in $\tau'$), we may use
Proposition~\ref{ThreeCasesCountableTightness}, Part A to find
finite sets $F_n^\gamma\subseteq D_n^\gamma$ in such a way that the
set $B=\bigcup\{F_n^\gamma:m<\omega\}$ is dense in $\tau'$. To
finish the construction, apply again part 1 of
Lemma~\ref{AngelosLemma2AboutMaximal} to get a topology
$\tau_{\gamma+1}$ which is the enlargement of $\tau'$ where $B$ is
dense open.

Let $\tau$ be the topology generated by $\bigcup\{\tau_\alpha:\alpha
<\mathfrak c\}$.   If the set $A$ is crowded in $\tau$ and
$A=A_\alpha $, then $A$ is also crowded in $\tau_\alpha  $. By
construction, $A=A_\alpha $ is open in $\tau_{\alpha +1}$ and so
even in $\tau$ (the second possibility in condition 4 cannot occur
because  $A$ cannot have isolated points  in $\tau_{\alpha +1}$. The
fact that every crowded subset of $\tau$ is open ensures that $\tau$
is a maximal topology \cite{vDMaximal}. If $\langle
D_n:n<\omega\rangle $ is a sequence of dense sets  in $\tau$, then
each $D_n$ is dense in each $\tau_\alpha $ and so there is some
$\beta<\mathfrak c$ such that each $D_n$ is dense open in
$\tau_\beta$. Since every $\omega$-sequence of subsets of $\omega$
is listed $\mathfrak c$-many times, there is an ordinal $\gamma\ge
\beta$ such that $\langle D_n:n<\omega\rangle=\langle
D_n^\gamma:n<\omega\rangle$. By condition 5 we get finite sets
$F_n^\gamma\subseteq D_n^\gamma=D_n$ in such a way that the set
$\bigcup\{F_n^\gamma:n<\gamma\}$ is dense open in $\tau_{\gamma+1}$
and so dense even in $\tau$. This shows that the space
$(\omega,\tau)$ is selectively separable.
\end{proof}

\vspace{3mm} Now we go back to D-separability. It turns out that, at
least in the countable case, maximal spaces are always D-separable.

\begin{theorem}\label{AngelosTheorem5AboutMax}
Let $X$ be submaximal. Then the following conditions are equivalent:
\begin{enumerate}
\item $X$ is D-separable;
\item $X$ is d-separable;
\item $X$ is $\sigma$-discrete;
\item $X$ is $\sigma$-closed discrete.
\end{enumerate}
\end{theorem}

The following is immediate:

\begin{corollary}\label{ImmediateFromAngelosFive}
A countable submaximal space is D-separable.
\end{corollary}

\begin{proof}[Proof of Theorem~\ref{AngelosTheorem5AboutMax}] First, we
prove the theorem for the special case when $X$ is crowded. The
implications (1)$\Rightarrow$(2) and
(4)$\Rightarrow$(3)$\Rightarrow$(2) are obvious.

(2)$\Rightarrow$(1)$\&$(4): Let $H=\bigcup_{n\in\omega} H_n$ be a
dense subspace of $X$ where each $H_n$ is discrete. Further, let
$(D_n : n\in\omega)$ be an arbitrary sequence of dense subsets of
$X$. For $n\in\omega$, put $D_n^\prime = D_n \cap (H\setminus H_n)$.
Ten $D_n^\prime$ is dense in $X$. Next, put $G_n=(D'_0\cap
\cdots\cap D'_{n})\setminus D'_{n+1}$. Each $G_n$ is closed discrete
being a subset  of the complement to the dense set $D'_{n+1}$.  By
construction, we have $G_n\subseteq D_n$ and
$D_0^\prime=\bigcup_{n\in\omega}G_n$.   Therefore,
$\bigcup_{n\in\omega}G_n$ is dense in $X$. Last, put $G_\omega = X
\setminus D_0^\prime$. Then $X=\bigcup_{n\leq\omega} G_n$ is a
countable union of closed discrete subspaces.

So we have proved the theorem for crowded $X$. It follows in
particular that (*) every countable crowded submaximal space is
D-separable. Now let $X$ be arbitrary submaximal space. Replace
every isolated point of $X$ with a copy of a countable crowded
regular maximal space. Call the resulting space $\tilde{X}$;
$\tilde{X}$ is crowded. It is easy to deduce from (*) that $X$ has
one of the properties (1) through (4) iff so does $\tilde{X}$. This
completes the proof.
\end{proof}

We will see that submaximal spaces can never be $D^+$-separable (see
Theorem $\ref{D+separableomegaresolvable}$). However, we don't know
the answer to the following question.

\begin{question}
{\rm Is there a countable submaximal space which is not
$DH$-separable?}
\end{question}

Arhangel'skii and Collins asked in \cite{ArhCollinsSubmaximal} if
all submaximal spaces are $\sigma$-discrete. Schr\"oder proved in
\cite{SchroderSubmaximal} that assuming V=L the answer is
affirmative. Thus we get:

\begin{corollary}
{\rm [V=L]} Every submaximal space is D-separable.
\end{corollary}

On the other hand, Kunen, Szymanski and Tall showed
\cite{KunenSzymanskiTall} that the existence of a measurable cardinal is consistent with ZFC iff the existence of a
Tychonoff crowded SIB space is consistent with ZFC. Further, Levy
and Porter proved the following:

\begin{proposition}\label{EqCondiditionsDueToLevyAndPorter}
{\rm (\cite{LevyPorter}, Proposition~3.1 and a remark after it)} The
following conditions are equivalent:\footnote{The equivalence of
conditions (6) and (7) is attributed in \cite{LevyPorter} to
Malykhin \cite{MalykhinResolvabilityProblemKatetov}.}
\begin{enumerate}
\item There exists a submaximal Hausdorff space which is not
$\sigma$-discrete;
\item There exists a crowded submaximal Hausdorff space which is not
$\sigma$-discrete;
\item There exists a maximal  space which is not
$\sigma$-discrete;
\item There exists a crowded submaximal Hausdorff space which is not
strongly $\sigma$-discrete;
\item There exists a maximal space which is not strongly
$\sigma$-discrete;
\item There exists a crowded SIB space;
\item There exists a crowded Hausdorff space $X$ such that every
real-valued function defined on $X$ is continuous at some point.
\end{enumerate}
\end{proposition}
It follows \cite{LevyPorter} that the existence of a crowded submaximal (or, equivalently, the existence of a
maximal) Hausdorff space which is not $\sigma$-discrete is equiconsistent with a measurable cardinal.  We see
from Theorem~\ref{AngelosTheorem5AboutMax} that the existence of a
submaximal space which is not D-separable can be added to the list
of conditions in Proposition~\ref{EqCondiditionsDueToLevyAndPorter}.
Therefore we get:

\begin{corollary}
The existence of a submaximal space which is not D-separable is
equiconsistent with a measurable cardinal.
\end{corollary}

As there exist (in ZFC) countable regular maximal non-selectively
separable spaces, Corollary~\ref{ImmediateFromAngelosFive} implies the existence of a space with some of the properties of Example $\ref{answerboaz}$.

\begin{corollary}\label{AngelosCorollary2AboutMax}
There exists a countable regular (maximal) D-separable
non-selectively separable space.
\end{corollary}

Unlike Example $\ref{answerboaz}$ such a space can never be sequential. Indeed, maximal spaces contain no non-trivial convergent sequences. However, Malykhin \cite{MalykhinMaximalGroups} has shown that maximal spaces can carry a group structure, so we wonder if Corollary $\ref{AngelosCorollary2AboutMax}$ can be improved in the following way:

\begin{question}
{\rm Is there a countable regular maximal non-selectively separable
topological group?}
\end{question}

Such an improvement can never be achieved in ZFC alone because Protasov \cite{ProtasovMaximalGroups} has shown that there are models of ZFC with no maximal topological groups.

We conclude by showing that a crowded submaximal space cannot be
D$^+$-separable. It was shown in \cite{BermanDowTopProc} that every
crowded SS$^+$ space is resolvable (and hence cannot be
submaximal). We note that the argument extends to D$^+$-separable
spaces.

\begin{theorem} \label{D+separableomegaresolvable}
Every crowded D$^+$-separable space is $\omega$-resolvable (and
hence non-submaximal).
\end{theorem}

\begin{proof}
Let $\sigma$ be the winning strategy for Two in the
game ${\sf G}_{\rm dis}({\mathcal D},{\mathcal D})$ on the space
$X$. It suffices to show that any dense $D\subset X$ contains two
disjoint dense subsets. Let $S_0=\sigma(D)$, $T_0=\sigma(D\setminus
S_0)$, $S_1=\sigma(D, D\setminus(S_0\cup T_0))$,
$T_1=\sigma(D,D\setminus(S_0\cup T_0\cup S_1))$,
$S_2=\sigma(D,D\setminus(S_0\cup T_0), D\setminus(S_0\cup T_0\cup
S_1\cup T_1))$, etc... Because $\sigma$ is winning the disjoint sets
$\bigcup_{n\in\omega}S_n$ and $\bigcup_{n\in\omega}T_n$ are dense.
\end{proof}

\begin{proposition}
If $X$ is D-separable and $\omega$-resolvable, then  there is a
family $\mathcal G$ of dense subspaces of $X$ such that $|{\mathcal
G}|=\mathfrak c$ and $G\cap G^\prime$ is nowhere dense, for every distinct
$G, G^\prime \in \mathcal G$.
\end{proposition}

\begin{proof}
Fix a family $\{Y_n : n\in\omega\}$ of pairwise
disjoint dense subsets of $X$. Also fix an almost disjoint family
$\mathcal A$ of infinite subsets of $\omega$ such that $|{\mathcal
A}|=\mathfrak c$. For each $A\in \mathcal A$ apply the definition of
D-separability to the family of dense subspaces $\{Y_n : n\in A\}$
to get discrete $F_{A,n} \subset Y_n$ such that $Z_A= \bigcup_{n\in
A} F_{A,n}$ is dense in $X$. It remains to note that whenever $A,
A^\prime \in \mathcal A$ are distinct, $Z_A \cap Z_{A^\prime}$ is
the union of finitely many discrete sets and thus nowhere dense.
\end{proof}

\begin{corollary}\label{CorollaryCountableExtraResolvable}
If $X$ is D-separable, $\omega$-resolvable and $\Delta(X)<\mathfrak c$
(in particular, if $X$ is countable), then $X$ is extra-resolvable.
\end{corollary}

\begin{corollary}
(1) If $X$ is  a crowded D$^+$-separable space then there is a family $\mathcal G$
of dense subspaces of $X$ such that $|{\mathcal G}|=\mathfrak c$ and
$G\cap G^\prime$ is nowhere dense, for every distinct $G, G^\prime \in
\mathcal G$.

(2) If $X$ is a crowded D$^+$-separable space and $\Delta(X)<\mathfrak c$ (in
particular, if $X$ is countable), then $X$ is extra-resolvable.
\end{corollary}

Recently, Garcia-Ferreira and Hru\v{s}ak \cite{GarciaHrushak}
constructed (within ZFC) a countable $\omega$-resolvable space which
is not extra-resolvable. It follows from
Corollary~\ref{CorollaryCountableExtraResolvable} that this provides
one more example of a countable space which is not D-separable.

\begin{example}
The Garcia-Ferreira-Hru\v{s}ak example of an $\omega$-resolvable non extra-resolvable space is an example of a countable non-$D$-separable space.
\end{example}


\vspace{3mm}{\bf Acknowledgment.} 
The authors express gratitude to Salvador Garcia-Ferreira, Paul
Gartside, Ronnie Levy and Boaz Tsaban for useful discussion.


\end{document}